\documentclass{article}
\usepackage{amsmath}
\usepackage{amsfonts}
\usepackage{longtable}
\usepackage{graphicx}

\newcommand{\JAA}{$J^{(3,1)}$}
\newcommand{\JAB}{$J^{(3,2)}$}
\newcommand{\JAC}{$J^{(3,3)}$}
\newcommand{\JAD}{$J^{(3,4)}$}
\newcommand{\JAE}{$J^{(3,5)}$}
\newcommand{\JAF}{$J^{(3,6)}$}
\newcommand{\JBA}{$J^{(4,5)}$}
\newcommand{\JBB}{$J^{(4,4)}$}
\newcommand{\JBC}{$J^{(4,3)}$}
\newcommand{\JBD}{$J^{(4,2)}$}
\newcommand{\JBE}{$J^{(4,1)}$}
\newcommand{\JBF}{$J^{(4,6)}$}
\newcommand{\JCA}{$J^{(5,2)}$}
\newcommand{\JCB}{$J^{(5,1)}$}
\newcommand{\JDA}{$J^{(6,1)}$}

\newcommand{\C}{\mathbb{C}}

\newcommand{\Z}{\mathbb{Z}}
\newcommand{\m}{\mathfrak{m}}

\newcommand{\Mat}{Mat}
\newcommand{\Sing}{Sing}

\newtheorem{thm}{Theorem}[section]
\newtheorem{lem}[thm]{Lemma}
\newtheorem{rem}[thm]{Remark}
\newtheorem{df}[thm]{Definition}
\numberwithin{equation}{section}

\newenvironment{proof}{\noindent
                       {\bf Proof:}}{
                       \begin{flushright}$\Box$\end{flushright}}
\newenvironment{entry}%
    {\begin{list}{}%
        {%
        \setlength{\labelwidth}{20pt}%
        \setlength{\labelsep}{5pt}%
        \setlength{\leftmargin}{25pt}%
        }%
    }%
{\end{list}}

\newcommand{\struttop}{\rule[0pt]{0pt}{12pt}}
\newcommand{\strutbot}{\rule[-4pt]{0pt}{4pt}}
\newcommand{\struta}{\rule[-6pt]{0pt}{18pt}}
\newcommand{\strutb}{\rule[-12pt]{0pt}{30pt}}

\begin{document}
\title{Simple Cohen-Macaulay Codimension 2 Singularities}

\author{\begin{tabular}{ll}
        Anne Fr\"uhbis-Kr\"uger & Alexander Neumer\\ 
        Inst. f. Alg. Geom. & Stiftung IHF\\ 
        Leibniz Univ. Hannover & Bremser Str. 79 \\ 
        30167 Hannover & 67063 Ludwigshafen/Rh.\\ 
        Germany & Germany \end{tabular}}


\maketitle

\begin{abstract}
In this article, we provide a complete list of simple isolated
Cohen-Macaulay codimension 2 singularities together with a list of
adjacencies which is complete in the case of fat point and space curve
singularities. 

\end{abstract}

\section{Introduction}

In this article, we determine a complete classification of simple 
Cohen-Macaulay codimension 2 singularities and in the case of
fat points and curves their complete adjacency list . 
Classification is done up to isomorphism of germs;
a singularity is called simple, if it can only deform into finitely
many different isomorphism classes -- in other words, if the modality
of the singularity is zero. Cohen-Macaulay codimension 2 singularities
are particularly important, as not all of them are complete intersections, but
they are, nevertheless, unobstructed and the theorem of Hilbert-Burch
provides a powerful tool for describing these singularities and their
deformations.  \\
\\
Arnold was the first to choose modality as the
criterion in his pioneering work \cite{Arn}. There he stated the famous
ADE-list of simple hypersurface singularities. About a decade
later Giusti \cite{Giu} gave a list of simple complete intersection
singularities; shortly thereafter Wall \cite{Wal} extended this to a
classification of unimodal isolated complete intersection
singularities. But there are singularities of modality 0 which are not 
complete intersections, as the following two examples show: the three
coordinate axes in $(\C^3,0)$, given by $xy,xz,yz$, can be checked to
be simple by direct computation, and the singularity in $(\C^6,0)$
given by the 2-minors of the $3 \times 2$ matrix, whose entries are
exactly the variables, does not permit any non-trivial deformations
($T^1 = 0$). In the case of space curve singularities, Giusti's list
was completed by the list of simple non-complete-intersection space
curves by Fr\"uhbis-Kr\"uger \cite{FK1} another decade later. For
readers' convenience, the above-mentioned lists of simple
singularities are included in the appendix, as we are computing
adjacencies to them.\\ 
\\
For the first task in this article, the classification, the main tools
are the Hilbert-Burch theorem (\cite{Bur}) and its generalization to 
deformations of singularities (cf. \cite{Sch}, \cite{Art}).
These allow us to describe the singularity by the presentation matrix
and the flat deformations by perturbations of this matrix. So the role,
which is played by the ideal in the case of complete intersection
singularities, is now taken over by the presentation matrix of this
ideal; this leads to reformulations of the $T^1$ and a finite
determinancy statement in terms of the presentation matrix as is shown
in detail in \cite{FK1}.\\
\\
In the proof of simplicity in each of the above-mentioned classifications
by Giusti and Fr\"uhbis-Kr\"uger, also a large number of adjacency
relations has been determined, but it was often unnecessary to consider all
adjacency relations; so the adjacency lists stated there cannot be
assumed to be complete. For instance, the adjacency list for the
simple space curves from Giusti's list was not
completed before the late 1990s (cf. \cite{S-V}), although it
had already been known to be incomplete for nearly a decade
(cf. \cite{Gor}), and even then the adjacencies to plane curve
singularities were not determined.\\ 
\\
The main tool for excluding adjacencies is semicontinuity of numerical
invariants, such as the Tjurina number, Milnor number and Delta
invariant. But even in the case of space curve singularities, this
tool is far from being powerful enough to decide all possible
cases. In this situation, determining the complete list of adjacencies
would not have been possible without systematic use of the computer
algebra system {\sc Singular} (\cite{Sin}) and, in particular, the
partial standard bases algorithm, a specialized tool for computing
simultaneously in families of singularities (cf. \cite{FK2},
\cite{FK3}). For each of the singularities, this allowed us to compute
the stratum in the base of the versal family, where the Tjurina number
is exactly one less than the one of the original singularity. Equipped
with the additional knowledge, which singularities of this particular
value of $\tau$ were appearing in the versal family, it was then
possible to reduce the cases, which had to be excluded by explicit
calculation, to just 3. In higher dimensions, the same method can be
applied to those singularities which are not part of a series, but this is not part of this article.\\
\\
With these lists, we hope to provide a set of examples of singularities
together with their adjacencies also in cases which are not hypersurfaces,
curves or surfaces. We should like to thank Gerhard Pfister and
the whole algebraic geometry group at the University of Kaiserslautern 
as well as the developers of the computer algebra system {\sc Singular} 
\cite{Sin} for many fruitful discussions. We would also like to thank
Jan Stevens for pointing out an omission in an earlier version of
this article and for several helpful remarks.

%

\section{Basics}
\label{Basics} For a detailed discussion of the methods used to study
germs of Cohen-Macaulay codimension 2 singularities with respect to
contact-equivalence, see \cite{FK1}. In short, we can say that using
the Hilbert-Burch theorem, all Cohen-Macaulay germs of 
codimension 2 can be expressed as the maximal minors of 
$(n+1)\times n$-matrices $M$ and vice versa. In the same way, flat 
deformations can be represented by perturbations of the matrix $M$ and 
any perturbation gives rise to a flat deformation (cf. \cite{Bur}, 
\cite{Sch}). 

Classification up to contact-equivalence means that two singularities
are considered equivalent, if the germs are isomorphic.\footnote{We
use the symbol $\sim_C$ to indicate contact-equivalence.} The action of
the contact-group translates directly to the application of coordinate
changes and row and column operations on $M$. A singularity is called
simple, if it can only deform into finitely many different equivalence
classes (types) of singularities.

\begin{df}[\cite{FK1}]
Let $M$ be a $(n+1) \times n$ matrix with entries in ${\mathbb C}\{x_1,
\dots,x_m\}$. $M$ is called quasihomogeneous of type $(D;a)\in
Mat((n+1)\times n;{\mathbb N}) \times {\mathbb N}^m$, if \\
a) all entries $M_{ij}$ are quasihomogeneous of degree $D_{ij}$ with
   respect to the weight vector a\\
b) there are relative row and column weights, i.e. \\
   \hspace{1cm} $D_{ij} - D_{ik} = D_{lj} - D_{lk} \;\; for\; all\; 1
   \leq i,l \leq n+1,\; 1 \leq j,k \leq n$ \\
Let $N$ be another $(n+1) \times n$ matrix with entries in ${\mathbb
  C}\{x_1,\dots,x_m\}$. The relative matrix weight of $N$ with respect
to $(D;a)$ is given by 
$$ v_{(D;a)}(N):=\inf_{j,i}\{v_a(N_{ij})-D_{ij}\} $$
\end{df}

\begin{lem}[\cite{FK1}]
Let $(X,0)$ be an isolated Cohen-Macaulay codimension 2 singularity
which is quasihomogeneous w.r.t. some weight vector $a$. Then it is 
possible to find a presentation matrix $M$ (describing the singularity
$X$) which is quasihomogeneous of type $(D;a)$ for a suitable $D \in
Mat((n+1)\times n;{\mathbb N})$. 
\end{lem}

For a consistent notation in the discussion, it is also necessary to
reformulate $T^1_{X,0}$ and the finite determinancy criterion in terms of
the presentation matrix\footnote{If we are considering a singularity $X,0$ 
in the notation of its presentation matrix $M$, we often also denote
$T^1_{X,0}$ by $T^1(M)$}: 

\begin{lem}[\cite{FK1}]
$T^1_{X,0}$ is given by
$$T^1_{X,0} \cong Mat(n+1,n;{\mathbb C}\{x_1,\dots,x_n\}) / (J(M) + Im(g))$$
where J(M) is the submodule generated by the matrices of the form
$$ \begin{pmatrix}
  \frac{\partial M_{11}}{\partial x_j} & \dots &
  \frac{\partial M_{1n}}{\partial x_j} \cr
  \vdots & & \vdots \cr
  \frac{\partial M_{(n+1)1}}{\partial x_j} & \dots &
  \frac{\partial M_{(n+1)n}}{\partial x_j} \end{pmatrix}\;\;\; \forall 1 \leq j
  \leq m$$
and $g$ is the map 
\begin{multline*}
\Mat(n+1,n+1;\C\{x_1,\ldots,x_m\})\oplus \Mat(n,n;\C\{x_1,\ldots,x_m\})
\\ \overset{g}{\rightarrow}  \Mat(n+1,n;\C\{x_1,\ldots,x_m\})
\end{multline*}
mapping $(A,B)  \mapsto  AM+MB$.
\end{lem}
 
By using the relative matrix weight, $T^1_{X,0}$ can be regarded as a
graded module $\bigoplus\limits_{v\in\Z}T^1_\nu(M)$. 

This can in turn be used to formulate an explicit determinancy
criterion for isolated quasihomogeneous Cohen-Macaulay codimension 2
singularities:  

\begin{lem}[\cite{FK1}]
Let $M$ be a $(n+1) \times n$ matrix with entries in the maximal ideal
of ${\mathbb C}\{x_1,\dots,x_m\}$, quasihomogeneous of type $(D;a)$
and defining an isolated singularity. Let $N$ be another $(n+1) \times
n$ matrix with entries in ${\mathbb C}\{x_1,\dots,x_m\}$, such that 
$$ v_{(D;a)}(N) > \beta = sup\{0,\alpha\} $$
where $\alpha = sup\{\nu \in {\mathbb Z} | T^1_{\nu}(M) \neq 0 \}$.\\
Then $M+N \sim_C M$.
\end{lem}

\section{Candidates in dimension $\geq 4$}
\label{geq4}
\subsection{Reduction of the problem}
\label{Reduction}
 Let $G_0$ be the $\C$-vector space of all quasihomogeneous $(n+1) \times n$
matrices of type $(D,a)$ for an arbitrary fixed positive integer $n$. 
For a generic matrix $M\in G_0$,  the kernel of the natural surjection
$G_0 \longrightarrow T^1_0(M)$ is generated by the set $S_1\cup
S_2\cup S_3$  where
\begin{eqnarray*}
S_1 & = & \left\{\left. b_{i,a_j}\frac{\partial M}{\partial x_j} \right| 1\leq i \leq r(a,a_j), 1\leq j \leq m\right\} \\
S_2 & = & \left\{\left. b_{i,D_{l1}-D_{j1}}Z_{jl}\right| 1\leq i \leq r(a,D_{l1}-D_{j1}), 1\leq j,l \leq n+1\right\} \\
S_3 & = & \left\{\left. b_{i,D_{1l}-D_{1j}}S_{jl}\right| 1\leq i \leq r(a,D_{1l}-D_{1j}), 1\leq j,l \leq n\right\}. \\
\end{eqnarray*}
Here $\{b_{1,d},\ldots,b_{r(a,d),d}\}$ is the set of monomials of weighted degree $d$, $r(a,d)$ its cardinality; $Z_{jl}$ denotes the $(n+1)\times n$-matrix,
having the $l$-th row of $M$ as its $j$-th row and all other entries $0$. In the same way, $S_{jl}$ is the matrix with the $l$-th column of $M$ as its $j$-th
column.

Recall that a singularity is called simple, if it can only deform into finitely many equivalence classes of singularities. Thus counting degrees of freedom shows
that a singularity defined by  an element $M\in G_0$ of type $(D;a)$ can only be simple if the dimension of $G_0$, which is just $\sum_{ij}r(a,D_{ij})$, does not
exceed the dimension of the above kernel, i.e. the number $s$ of linearily independent elements of $S_1,S_2$ and $S_3$. Since we always have two relations, the
Euler relation and $\sum S_{ii}=\sum Z_{jj}$, this means that the inequality $\#S_1+\#S_2+\#S_3-2\geq \sum_{ij}r(a,D_{ij})$ has to hold.

\begin{lem}
A $(n+1)\times n$-matrix $M$ can only define a simple isolated
codimension 2 singularity in $(\C^m,0)$ if $n<3$ and $m<7$.
\end{lem}
\begin{proof}
We will use a modified version of the above counting argument to prove
that the number of different variables occuring in the 1-jet of $M$
has to be greater than $n^2+n-2$ whereas the total number of variables
may not exceed $n^2+n$.

{\em Step 1 (An upper bound for the number of variables)}
Let us first suppose that the total number of variables $m$ exceeds
$n^2+n$. Consider a generic $(n+1)\times n$-matrix, quasihomogeneous
w.r.t. the weights $a=(1,\ldots,1)$ and the degrees $D_{ij}=1$. The
matrix contains $n^2+n$ different linear entries and is hence
equivalent to 
$$\begin{pmatrix} x_1 & \dots & x_n \\
           \vdots & & \vdots \\
            x_{n^2+1} & \dots & x_{n(n+1)}
\end{pmatrix}.$$
This matrix does not allow any non-trivial perturbations and any other
matrix of the corresponding size in $m$ variables is adjacent to it;
on the other hand, a direct computation shows that it defines a
non-isolated singularity. Therefore isolated singularities can
only occur for $m \leq n^2+n$.

{\em Step 2 (A counting argument on the 1-jet)} Let $N$ be a generic
$(n+1)\times n$-matrix, quasihomogeneous w.r.t. the weights
$a=(1,\ldots,1)$ and the degrees 
$D_{ij}=1$. For these weights, the kernel of the map from $G_0$ to $T^1_0(N)$ in the above argument is generated by $\#S_1+\#S_2+\#S_3-2=m^2+(n+1)^2+n^2-2$
elements. Comparing this to $m(n^2+n)$, the dimension of $T^1_0(N)$, we get the inequality $m^2+(n+1)^2+n^2-2\geq m(n^2+n)$ as a necessary condition for $N$ being
simple. This simplifies to $m^2>(m-2)(n^2+n)$. First of all, we see directly from this inequality that for $n>2$, $m$ has to be at least $10$. Using this
additional information, we can now simplify our condition to
\begin{multline*}
m+3>m+2+\frac{4}{m-2}=\frac{(m+2)(m-2)+4}{m-2}=\\=\frac{m^2}{m-2}>\frac{(m-2)(n^2+n)}{m-2}=n^2+n.
\end{multline*}
Since $N$ was a generic matrix of this type, the 1-jet of any matrix $M$ of appropriate size with entries in the maximal ideal of $\C\{x_1,\ldots,x_m\}$ is
adjacent to $N$.

Combining this with the result of step 1, we see that $n^2+n-2 \leq m
\leq n^2+n$, but we still do not have a bound for the number of
variables actually appearing in the 1-jet nor do we have a bound for
the size of the matrix. To obtain these two, we now pass to another
set of weights: Let $N$ be a generic
$(n+1)\times n$-matrix, quasihomogeneous w.r.t. the
weights\footnote{Any matrix whose 1-jet only involves $p$ of the $m$
variables is adjacent to this matrix and hence this is the set of
weights to consider for determining the least number of variables
appearing in the 1-jet of a simple singularity of given size $n \times
(n+1)$ and number of variables $m$.}
$$a=(\underbrace{2,\dots,2}_{p},\underbrace{1,\dots,1}_{m-p}) \;\;\;\;
  D_{ij}=2.$$  
For these weights, we obtain
$$\#S_1+\#S_2+\#S_3-2=n^2+(n+1)^2+p(p+(m-p) \cdot \frac{m-p+1}{2})+(m-p)^2$$
and
$${\rm dim}_{\mathbb C} T^1_0(N)=n(n+1)(p+(m-p) \cdot \frac{m-p+1}{2}).$$
Plugging in each of the three possible values of $m$, we obtain the
following table whose entries are the difference of the number of
degrees of freedom and the number of possible entries minus 2. Hence
simple singularities can only occur if the entry has a positive value.

\begin{tabular}{|c|c|c|c|}
\hline
 & $m=n^2+n-2$ & $m=n^2+n-1$ & $m=n^2+n$ \\
\hline
$p=n^2+n$ & - & - & $2n^2+2n-1$\\
\hline
$p=n^2+n-1$ & - & $n^2+n$ & $n^2+n$\\
\hline
$p=n^2+n-2$ & 3 & 2 & 1\\
\hline
$p=n^2+n-3$ & $-n^2-n+6$ & $-n^2-n+3$ & -$n^2-n-1$\\
\hline
$p=n^2+n-3$ & $-2n^2-2n+7$ & - & -\\
\hline
\end{tabular}\\

We immediately see that at least $n^2+n-2$  variables have to appear
in the 1-jet resp. $n^2+n-3$ in the case $n=2$, $m=4$.
 
{\em Step 3 (Excluding non-isolated singularities)} The conditions
obtained in steps 1 and 2 imply that (for $n>2$) $j_1M$ has to contain
at least $n^2+n-2$ different variables. But then $j_1M$ is
contact-equivalent to a matrix of the form 
$$j_1M \sim_C \begin{pmatrix}
x_{11} & \ldots & x_{(n-1)1} & x_{n1} & x_{(n+1)1} \\
\vdots & & & & \vdots \\
x_{1(n-1)} & \ldots & x_{(n-1)(n-1)} & \alpha & x_{(n+1)(n-1)} \\
x_{1n} & \ldots & x_{(n-1)n} & x_{nn} & \beta
\end{pmatrix}
$$
with $\alpha, \beta\in\m$, and by another coordinate change we see
that $M$ is of the same form. 

Direct computation shows that the ideal of the singular locus of $M$
is contained in
$\left<x_{11},\ldots,x_{1n},x_{21},\ldots,x_{2n}\right>$. This
component is obviously not zero-dimensional, and hence the singularity
defined by $M$ cannot be isolated. 
\end{proof}

We can now restrict our consideration to the case of $(n+1)\times
n$-matrizes with $n\leq 2$. In the case $n=1$, the singularity is a
complete intersection; a complete classification of simple isolated
singularities in this situation can be found in \cite{Giu}. 

In the remaining case of $3\times 2$-matrizes with $n=2$, the
calculations of step 2 of the preceding proof imply that
simple isolated singularities can only occur for:
\begin{center}
\begin{tabular}{c|c}
dimension $m$ & possible 1-jet candidates \strutbot \\
\hline \hline
 3 & 3 variables \struttop \\
 4 & 3, 4 variables \\
 5 & 4, 5 variables \\
 6 & 4, 5, 6 variables \\
\end{tabular}
\end{center}

\paragraph{Remark:} Unfortunately, this method gives no bound for the
case of fat points in $(\C^2,0)$. Thus we will study fat points
seperately in section \ref{FatPoints}. 

As the next step, we will classify the possible candidates of 1-jets with 4 or more variables:

\subsection{1-jet candidates}
First we will classify the possible candidates of 1-jets with 4 or
more variables:  
\begin{lem}
\label{jetlemma} Let $M$ be a $3\times 2$-matrix with entries in the maximal ideal of $\C\{x_1,\ldots,x_m\}$. Then $j_1M$ is contact-equivalent to one of the jets
in the following tables

\vspace{0,5cm}

\begin{center}
\parbox[top][7.2cm][t]{4cm}{

\vspace{0,5cm}

\begin{tabular}{|c|c|}
\hline
\multicolumn{2}{|c|}{6 variables\struta} \\
\hline \strutb\JDA & $\begin{pmatrix}
    x & y & v \\
    z & w & u
    \end{pmatrix}$ \\
\hline
\end{tabular} \\

\vspace{1cm}

\begin{tabular}{|c|c|}
\hline
\multicolumn{2}{|c|}{5 variables\struta} \\
\hline \strutb \JCB & $\begin{pmatrix}
    x & y & v \\
    z & w & x
    \end{pmatrix}$ \\
\hline \strutb \JCA & $\begin{pmatrix}
    x & y & v \\
    z & w & 0
    \end{pmatrix}$ \\
\hline
\end{tabular}}
\parbox[top][7.2cm][t]{4cm}{
\begin{tabular}{|c|c|} \hline
\multicolumn{2}{|c|}{4 variables\struta} \\
\hline \strutb \JBE & $\begin{pmatrix}
    w & y & x \\
    z & w & y
    \end{pmatrix}$ \\
\hline \strutb \JBD & $\begin{pmatrix}
    w & y & x \\
    z & w & 0
    \end{pmatrix}$ \\
\hline \strutb \JBC & $\begin{pmatrix}
    0 & y & x \\
    z & w & 0
    \end{pmatrix}$ \\
\hline \strutb \JBB & $\begin{pmatrix}
    x & y & z \\
    z & w & 0
    \end{pmatrix}$ \\
\hline \strutb \JBA & $\begin{pmatrix}
    x & y & 0 \\
    z & w & 0
    \end{pmatrix}$ \\
\hline \strutb \JBF & $\begin{pmatrix}
    x & y & z \\
    w & 0 & 0
    \end{pmatrix}$ \\
\hline
\end{tabular} }
\end{center}
or is contact-equivalent to a 1-jet containing only 3 or less variables.
\end{lem}

\vspace{10pt}

\noindent To simplify notation, we will abbreviate the case $k$ in
dimension $d$ as $J^{(d,k)}$. 

\vspace{0,5cm}

\begin{proof}
Because the arguments for all these cases work in a similar way, we will concentrate on the presentation of the case where $j_1M$ contains exactly $5$ variables.

By applying several coordinate changes and row and column operations we may assume w.l.o.g. that
$$j_1M\sim_C\begin{pmatrix}x & y & v \\ z & w & \alpha\end{pmatrix}\quad \text{ with $\alpha\in\m$.}$$
If $\alpha=0$, we have a 1-jet of type \JCA. Otherwise we can write $\alpha=\alpha_1x+\alpha_2y+\alpha_3z+\alpha_4w+\alpha_5v$ with $\alpha_i\in\C$. By
contact-equivalence, we get
$$j_1M\sim_C\begin{pmatrix}x & y & v \\ z & w & \alpha'_1x+\alpha'_2y\end{pmatrix}\quad\text{ with $\alpha'_1,\alpha'_2\in\C$.}$$
We may now assume $\alpha'_1\neq0$ (by exchanging the first and second column and exchanging the roles of $x$ and $y$ respectively of $z$ and $w$) and obtain
$$j_1M\sim_C\begin{pmatrix}x & y & v \\ z & w & x\end{pmatrix}.$$
\end{proof}

The possible 1-jets which contain exactly three variables were already
determined in \cite{FK1}:\\

\begin{center}
\begin{tabular}{|c|c||c|c|} \hline
\multicolumn{4}{|c|}{{\it 3 variables}\struta} \\
\hline \strutb \JAA & $\begin{pmatrix}
    z & y & x \\
    0 & x & y
    \end{pmatrix}$ &
 \strutb \JAB & $\begin{pmatrix}
    z & y & 0 \\
    0 & x & y
    \end{pmatrix}$\\
\hline \strutb \JAC & $\begin{pmatrix}
    z & y & 0 \\
    y & x & z
    \end{pmatrix}$ &
 \strutb \JAD & $\begin{pmatrix}
    z & 0 & 0 \\
    0 & x & y
    \end{pmatrix}$\\
\hline \strutb \JAE & $\begin{pmatrix}
    z & y & x \\
    y & 0 & 0
    \end{pmatrix}$ &
 \strutb \JAF & $\begin{pmatrix}
    z & 0 & x \\
    0 & z & y
    \end{pmatrix}$\\
\hline
\end{tabular}
\end{center}

\vspace{10pt}

The preceeding two lemmata provide a classification of all 1-jets that may
occur in simple isolated Cohen-Macaulay codimension 2 singularities in
$(\C^p,0)$, $p\geq 4$. Now we will check whether these 1-jets lead to
simple singularities. Since we already know that we will not get
simple singularities for all 1-jets in every dimension, we will start
by regarding the candidates in the smallest dimension. 

\subsection{Singularities in $({\mathbb C}^4,0)$}
Because we have seen that 1-jets containing only 2 or fewer variables
cannot be simple in dimension 4, we have to consider only the 1-jets
with 3 and 4 variables. 

\begin{thm}
The following table shows the list of simple isolated Cohen-Macaulay
codimension 2 singularities in $(\C^4,0)$.
{\small{
\begin{center}
\begin{longtable}{|c|c|l|l|c|l|}
\hline
    Jet-Type
    & Type
    & Presentation Matrix
    &
    & $\tau$ 
    & Name of Triple
\\  &
    &
    &
    &
    & Point in \cite{Tju}
\\ \hline
\JBE
    & $\Lambda_{1,1}$
    & $\begin{pmatrix}w & y & x \\ z & w & y \end{pmatrix}$
    &
    & 2
    & $A_{0,0,0}$
\\ \hline
\JBD
    & $\Lambda_{k,1}$
    & $\begin{pmatrix}w & y & x  \\ z & w & y^k \end{pmatrix}$
    & $k\geq 2$
    & $k+1$
    & $A_{0,0,k-1}$
\\ \hline
\JBC
    & $\Lambda_{k,l}$
    & $\begin{pmatrix}w^l & y & x \\ z & w & y^k \end{pmatrix}$
    & $k\geq l\geq 2$
    & $k+l$
    & $A_{0,l-1,k-1}$
\\ \hline
\JBB
    &
    & $\begin{pmatrix}z & y & x \\ x & w & y^2+z^k \end{pmatrix}$
    & $k\geq 2$
    & $k+3$
    & $C_{k+1,0}$
\\ \hline
    &
    & $\begin{pmatrix}z & y & x \\ x & w & yz+y^kw \end{pmatrix}$
    & $k\geq 1$
    & $2k+4$
    & $B_{2k+2,0}$
\\
    &
    & $\begin{pmatrix}z & y & x \\ x & w & yz+y^k \end{pmatrix}$
    & $k\geq 3$
    & $2k+1$
    & $B_{2k-1,0}$
\\ \hline
    &
    & $\begin{pmatrix}z & y & x \\ x & w & z^2+yw \end{pmatrix}$
    &
    & 7
    & $D_0$
\\ \hline
    &
    & $\begin{pmatrix}z & y & x \\ x & w & z^2+y^3 \end{pmatrix}$
    &
    & 8
    & $F_0$
\\ \hline
\JAA
    &
    & $\begin{pmatrix} z & y+w^l & w^m \\ w^k & y & x \end{pmatrix}$
    & $k,l,m \geq 2$
    & $k+l+m-1$
    & $A_{k-1,l-1,m-1}$
\\ \hline
\JAB
    &
    & $\begin{pmatrix}z & y & x^l+w^2 \\ w^k & x & y\end{pmatrix}$
    & $k,l\geq 2$
    & $k+l+2$
    & $C_{l+1,k-1}$
\\ \hline
    &
    & $\begin{pmatrix}z & y+w^l & xw \\ w^k & x & y \end{pmatrix}$
    & $k,l\geq 2$
    & $k+2l+1$
    & $B_{2l,k-1}$
\\
    &
    & $\begin{pmatrix}z & y & xw+w^l \\ w^k & x & y \end{pmatrix}$
    & $k\geq 2, l\geq 3$
    & $k+2l$
    & $B_{2l+1,k-1}$
\\ \hline
    &
    & $\begin{pmatrix}z & y+w^2 & x^2 \\ w^k & x & y \end{pmatrix}$
    & $k\geq 2$
    & $k+6$
    & $D_{k-1}$
\\ \hline
    &
    & $\begin{pmatrix}z & y & x^2+w^3 \\ w^k & x & y \end{pmatrix}$
    & $k\geq 2$
    & $k+7$
    & $F_{k-1}$
\\ \hline
\JAC
    &
    & $\begin{pmatrix}z & y & xw+w^k  \\ y & x & z \end{pmatrix}$
    &
    & $3k+1$
    & $H_{3k}$
\\
    &
    & $\begin{pmatrix}z & y & xw \\ y & x & z+w^k \end{pmatrix}$
    &
    & $3k+2$
    & $H_{3k+1}$
\\    &
    & $\begin{pmatrix}z & y & xw \\ y+w^k & x & z \end{pmatrix}$
    &
    & $3k+3$
    & $H_{3k+2}$
\\ \hline
    &
    & $\begin{pmatrix}z & y & w^2  \\ y & x & z+x^2 \end{pmatrix}$
    &
    & 8
    & 
\\ \hline
    &
    & $\begin{pmatrix}z & y & x^3+w^2  \\ y & x & z \end{pmatrix}$
    &
    & 9
    &
\\ \hline
    &
    & $\begin{pmatrix}z & y & x^2  \\ y & x & z+w^2 \end{pmatrix}$
    &
    & 9
    &
\\ \hline
\end{longtable}
\end{center}}}
\end{thm}

\begin{proof}
We will consider each of the possible 1-jets from lemma \ref{jetlemma}
seperately. As the proof that a singularity cannot be simple by
variants of the counting argument always has the same structure, we
only list the non-simple singularities and the respective weights here:
\begin{center}
\begin{longtable}{|c|l|c|c|c|}
\hline
    Jet-Type
    & Presentation Matrix
    & $\underline{a}$
    & $D$
    & $\tau$
\\ \hline
    \JBB
    & $\begin{pmatrix}z & y & x \\ x & w & z^2+y^4 \end{pmatrix}$
    & $\begin{pmatrix}2 & 1 & 3 \\ 3 & 2 & 4 \end{pmatrix}$
    & $\begin{pmatrix}3 & 1 & 2 & 2 \end{pmatrix}$
    & 11
\\ \hline
    & $\begin{pmatrix}z & y & x \\ x & w & y^3+z^3 \end{pmatrix}$
    & $\begin{pmatrix}1 & 1 & 2 \\ 2 & 2 & 3\end{pmatrix}$
    & $\begin{pmatrix}2 & 1 & 1 & 2 \end{pmatrix}$
    & 10
\\ \hline
    \JBA
    & $\begin{pmatrix}z & y & \alpha \\ x & w & \beta \end{pmatrix}$
    &    $\begin{pmatrix}1 & 1 & 2 \\ 1 & 1 & 2\end{pmatrix}$
    &    $\begin{pmatrix} 1 & 1 & 1 & 1 \end{pmatrix}$
    & 13
\\ \hline
    \JBF
     & $\begin{pmatrix} x & y & z \\ w & \alpha & \beta \end{pmatrix}$
     & $\begin{pmatrix}1 & 1 & 1 \\ 2 & 2 & 2 \end{pmatrix}$
     & $\begin{pmatrix} 1 & 1 & 1 & 2 \end{pmatrix}$
     & 9
\\ \hline
    \JAB
    & $\begin{pmatrix}z & y & \alpha \\ w^2 & x & y+w^2 \end{pmatrix}$
    & $\begin{pmatrix}3 & 2 & 3 \\ 2 & 1 & 2 \end{pmatrix}$
    & $\begin{pmatrix}1 & 2 & 3 & 1\end{pmatrix}$
    & 11
\\ \hline
    & $\begin{pmatrix}z & y+w^3 & x^2+\alpha w^4 \\ w^2 & x & y \end{pmatrix}$
    & $\begin{pmatrix} 3 & 3 & 4 \\ 2 & 2 & 3 \end{pmatrix}$
    & $\begin{pmatrix} 2 & 3 & 3 & 1 \end{pmatrix}$
    & 12
\\ \hline
    \JAD
    & $\begin{pmatrix} y & z & \alpha \\ \beta & y+\gamma & z+\delta\end{pmatrix}$
     & $\begin{pmatrix} 2 & 2 & 2 \\ 2 & 2 & 2 \end{pmatrix}$
     & $\begin{pmatrix} 1 & 2 & 2 & 1 \end{pmatrix}$
     & 15
\\ \hline
    \JAF
    & $\begin{pmatrix} z & y & x \\ \alpha & \beta & y+\gamma \end{pmatrix}$
     & $\begin{pmatrix} 4 & 6 & 4 \\ 6 & 8 & 6 \end{pmatrix}$
     & $\begin{pmatrix} 4 & 6 & 4 & 3 \end{pmatrix}$
     & 15
\\ \hline
\end{longtable}
\end{center}

\begin{entry}
\item[\JBE:] Let 
$$M=\begin{pmatrix} x & y & z \\ z & w & y \end{pmatrix}.$$ 
By direct computation, $T^1(M)$ is generated by
$$\begin{pmatrix}0&0&1\\0&0&0\end{pmatrix}\;\; {\rm and} 
  \begin{pmatrix}0&0&0\\0&0&1\end{pmatrix}.$$ 
Thus any deformation of $M$ is of the form 
$$M'(\alpha,\beta)=\begin{pmatrix}x&y&z+\alpha\\z&w&y+\beta\end{pmatrix}$$
with $\alpha,\beta\in\C$. By direct computation, the singular locus of
$M'(\alpha,\beta)$ is 
\begin{multline*}
\Sing(M'(\alpha,\beta)) = \left<\beta ^2,\alpha \beta ,\alpha ^2,w\beta ,5w\alpha -\beta ^2,w^2,5z\beta +4\alpha \beta ,\right. \\
5z\alpha +3\alpha ^2,zw+w\alpha ,2z^2+3z\alpha +\alpha ^2,2y\beta -zw+\beta ^2,3y\alpha +2z\beta +4\alpha \beta , \\
yw+w\beta , yz+y\alpha +z\beta +\alpha \beta ,2y^2+3y\beta +\beta ^2,x\beta -5z^2-3z\alpha,\\
\left. x\alpha ,xw-4y-2y\alpha -2z\beta -\alpha \beta ,xz+x\alpha ,xy+2z^2+z\alpha ,x^2 \right>
\end{multline*}
Because this ideal contains $\alpha^2$ and $\beta^2$,
$M'(\alpha,\beta)$ defines a smooth surface if either $\alpha$ or
$\beta\neq0$. Hence $M$ defines a simple singularity we will call
$\Lambda_{1,1}$. Since $M$ is 1-determined, every matrix of type \JBE is
contact-equivalent to $M$. 

\item[\JBD:] In this case $M$ is contact-equivalent to a matrix
$$\begin{pmatrix}w & y & x \\ z & w + \alpha & \beta \end{pmatrix}$$
with $\alpha,\beta\in\m^2$. Moreover, we can get rid of $\alpha$  and
all terms of $\beta$ involving $x,z,w$ by row and column operations
and by appropriate coordinate changes of $w$. This leads to the matrix
$$\begin{pmatrix}w & y & x \\ z & w & y^k\cdot(1+\gamma(y))
  \end{pmatrix}$$
where $\gamma \in {\mathfrak m}$. Dividing the last column by the unit
$(1+\gamma(y))$ and performing an appropriate coordinate change in
$x$, we then obtain the desired structure of the matrix.\\
For the proof of simplicity, we can proceed in the same way as in the
previous case and obtain that only adjacencies to singularities of the
same series with lower $\tau$, of the type \JBE and of the $A$-series
(i.e. $(w,y^k+xz)$) can appear.

\item[\JBC:] This case is strictly analogous to the case \JBD, with
the only difference that we obtain at most adjacencies to the
singularities of the same series with lower $\tau$, to the series \JBD, to
the singularity \JBE and to the A-series.

\item[\JBB:] A matrix with 1-jet of type \JBB is of the structure
$$\begin{pmatrix} z & y & x+\alpha \\ x & w & \beta \end{pmatrix}$$
where $\alpha,\beta \in {\mathfrak m}^2$. All terms of $\alpha$ can be
cancelled in the same way as in the case \JBD. Regarding $\beta$ only terms
in $y$ and $z$ and terms of the form $y^kw$ cannot be killed. By the
table of non-simple singularities, we can conclude that $\beta$ cannot
be of order 3 or higher. A direct calculation then shows that the
2-jet of $\beta$ has to be one of the following 7: $y^2+z^2$, $yz+yw$,
$z^2+yw$, $y^2$, $yz$, $z^2$ and $yw$. In the first three cases the
corresponding 2-jet of the matrix is already 2-determined which
implies that each gives rise to exactly one singularity. In the fourth
case, the only monomial of higher degree which may occur is a power of
$z$, leading to the first series in the list. In the 5th case, pure
powers of $y$ and terms $y^kw$ are the only
terms that cannot be cancelled, but for determinacy reasons more than
one of them cannot occur simultaneously; this gives rise to the second
series. In the last two cases, the only possibility which is not
excluded by the list of non-simple singularities (case \JBB, lines 1
and 2) is $z^2+y^3$.\\
For the proof of simplicity, we need to study two different questions
here: first of all, we have to find out whether some singularities
from the series are adjacent to non-simple ones and secondly, we have
to find out whether any of these singularities can deform into a
non-simple one of a different 1-jet. The first question can be
answered quite easily by observing that whenever a term $y^2$ is
present, the singularity is in the first series, and whenever a term
$yz$, but no $y^2$ is present, the term $z^2$ may be killed by a
coordinate transform in $y$ and a subsequent column operation on the
second column.\\
The second question involves a simple, but rather lengthy explicit
calculation which shows that a singularity of type \JBB can only
deform into singularities of types $A_k$, $D_k$, \JBE, \JBD and \JBC:
More precisely, we consider the versal family, deduce conditions
deciding when a point in the base space allows a singularity and then
determine the occurring kinds of singularities. The versal family in
our case is 
$$\begin{pmatrix} z & y & x +\alpha \\ x & w & \beta + \gamma y
+\delta w + \varepsilon z + p(y,z,w)\end{pmatrix},$$
$\alpha,\dots,\varepsilon \in \C$ and $p \in {\mathfrak m}^2$. As the
whole calculation is rather lengthy we 
only sketch a part of it, namely the case $\alpha = 0$: For fixed
$\beta, \dots, \varepsilon$, singularities may only occur  at points
where the order of the lower right hand entry of the matrix is at
least 1. By the structure of the matrix, we see that for $\beta
\neq 0$ at most a $D_k$ singularity ($x-zw,zw^2-\dots$) may occur, if
$y$ is non-zero at this point, at most an $A_k$ singularity
($w-yx,x^2+y^2+\dots$ or $w-yx,x^2-yz$) if the $z$ coordinate is
non-zero, and no singularities at other points. If $\beta=0$ and
the 1-jet is not of type \JBB, then $\gamma \neq 0$ implies that we are
dealing with a singularity of type \JBE, $\gamma=0,\delta \neq 0$
leads to \JBD and $\gamma=\delta=0,\varepsilon \neq 0$ leads to \JBC.

\item[\JBA:] As the generic matrix of this jet appears among the
non-simple singularities, we cannot get any candidates here.

\item[\JBF:] No simple singularities possible, same argument as for \JBA.

\item[\JAA:] In this case, the matrix is contact-equivalent to
$$\begin{pmatrix}z & y+\alpha & x+\beta \\ \gamma & x & y
  \end{pmatrix}$$ 
where $\alpha,\beta,\gamma \in {\mathfrak m}^2$; moreover, we can
achieve that $\alpha$, $\beta$ and $\gamma$ only contain terms in $w$
by appropriate row and column operations and coordinate changes in
$x$, $y$ and $z$. Therefore the matrix can be written as
$$\begin{pmatrix}z & y+a_1w^{k_1} & x+a_2w^{k_2} \\ a_3w^{k_3} & x & y
  \end{pmatrix}$$
where $a_1,a_2,a_3 \in \C\{w\}$ and either $a_i=0$ or $a_i(0)\neq 0$ for
$1 \leq i \leq 3$. Here it is most convenient to pass to an equivalent
way of writing the one jet (by a coordinate change $x\mapsto \frac{1}{2}(x+y)$ 
and $y \mapsto \frac{1}{2}(x-y)$ followed by appropriate row and column 
operations), in which the matrix can then be stated as
$$\begin{pmatrix}z & y+b_1w^{k_1} & b_2w^{k_2} \\ b_3w^{k_3} & y & x
  \end{pmatrix}$$
As the matrix would not describe a finitely determined singularity if 
any of the three $b_i$ were zero, we easily achieve that $b_2=b_3=1$
and with a little more work that also $b_1=1$. Writing the matrix in this
way, it is obviously quasihomogeneous and we hence write this normal
form in the list. On the other hand, transforming it back into the other form
facilitates comparisons in subsequent adjacency calculations:
$$\begin{pmatrix} z & y+w^l+w^m & x+w^l-w^m \\ w^k & x & y \end{pmatrix}$$
To prove simplicity, we proceed in the same way as outlined in \JBB
and obtain that the singularity can at most be adjacent to the
following singularities/series of singularities: $A$-series, \JBE,
\JBD and \JBC. Since all members of these series are simple as we
already proved, these singularities are simple, as well.

\item[\JAB:] By the same argument as in the case \JAA, the matrix has
to be of the structure
$$\begin{pmatrix}z & y +\alpha & \beta \\ w^k & x & y\end{pmatrix}$$
where $\alpha,\beta \in {\mathfrak m}^2$, $\alpha$ only involving
terms in $w$, $\beta$ only involving terms in $x$ and $w$. A matrix of
this structure does not define an isolated singularity if neither
$\alpha$ nor $\beta$ contain a pure power in $w$. Moreover, if the
order of $\beta$ is at least 3, the singularity cannot be simple due
to an adjacency to the non-simple singularities (\JAB line 1). Let 
us start with the case that $\beta$ 
contains the term $w^2$, which may w.l.o.g. be written as:
$$\begin{pmatrix}z & y +\alpha & a_1x^{k_1}+a_2x^{k_2}w+w^2\\
                 w^k & x & y\end{pmatrix}$$
where $a_1,a_2 \in \C\{x,w\}$ and either $a_i=0$ or $a_i(0)\neq 0$ for
each of them. By an appropriate coordinate change in $w$ followed by
substracting an appropriate multiple of the 2nd column from the first
one and another coordinate change in $z$, we can get rid of the term
$a_2x^{k_2}w$ changing of course the term $a_1x^{k_1}$ to some
$\tilde{a_1}x^{\tilde{k_1}}$. By determinacy, the terms of $\alpha$
may be omitted which leads to the first series.\\
If, on the other hand, the term $w^2$ is not present in
$\beta$, we are dealing with a matrix of the structure
$$\begin{pmatrix}z & y + b_1w^{l_1} & a_1x^2+a_2xw+b_2w^{l_2} \\
                 w^k & x & y\end{pmatrix}$$
where $a_1 \in \C\{x,w\}$, $a_2,b_1,b_2 \in \C\{w\}$ and at least one
of $a_1(0)$ and $a_2(0)$ and one of $b_1(0)$ and $b_2(0)$ non-zero. If
$a_2(0)\neq 0$, we can cancel the other term by an appropriate
coordinate change in $w$ (and, of course,  cleaning up as before) and
obtain the second series, because the term of higher relative matrix
weight out of $b_1w^{l_1}$ and $b_2w^{l_2}$ can be killed by
determinacy. If $a_2(0)=0$, $l_1=2$ and $b_1(0) \neq 0$, we can cancel
$b_2w^{l_2}$ again due to determinacy and obtain the third series.
If $a_2(0)=0$, no monomial $w^2$ appears in $b_1w^{l_1}$, $l_2=3$ and
$b_2(0)\neq 0$, then we obtain the 4th series. Otherwise, i.e. if
$a_2(0)=0$, no $w^2$ term appears in $b_1w^{l_1}$ and no $w^3$ term in
$b_2w^{l_2}$, then the singularity cannot be simple as 
it is adjacent to the non-simple singularity (\JAB, line 2).\\
The proof of simplicity can be done as in \JBB leading to the
following possible adjacencies: $D$-series, $A$-series, \JAA (with
$w^2$ in the upper right hand entry), \JBE, \JBD,
\JBC and \JBB. As not all singularities of series \JBB are simple, we
need to consider these adjacencies more closely: A simple explicit
calculation shows that these adjacencies, which are obtained by perturbing
with $w$ in the lower left-hand entry, only allow adjacencies to the
first series of \JBB for singularities of the first series,  
to the very first singularity in the same series for the second series
and to the singularity of Tjurina number 7 resp. 8 (\JBB, line 4
resp. 5) for the last remaining series. All the above mentioned
singularities of type \JBB are simple which in turn implies that the 3
series are simple as well. 

\item[\JAC:] A matrix with this 1-jet is of the structure
$$\begin{pmatrix}y & z & \alpha\\ 
                 x & y+\beta & z+\gamma\end{pmatrix}$$
where $\alpha,\beta,\gamma \in {\mathfrak m}^2$. By the table of
non-simple singularities, we know from the table of non-simple
singularities (\JBB,line 1) that the
generic matrix with the weights  $\underline{a}=(1,2,3,2)$ and 
$$D=\begin{pmatrix}2&3&4\\1&2&3\end{pmatrix}$$ cannot be simple. This
implies, that simple singularities can only occur in the following 6
cases  
\begin{enumerate}
\item[(a)] $\alpha=xw+w^k$
\item[(b)] $\alpha=xw$, $\gamma=w^k$
\item[(c)] $\alpha=xw$, $\beta=w^k$
\item[(d)] $\alpha=w^2$, $\gamma=x^2$
\item[(e)] $\alpha=x^3+w^2$
\item[(f)] $\alpha=x^2$, $\gamma=w^2$
\end{enumerate}
It can be seen directly by the usual determinacy argument that,
whenever $\beta$ resp. $\gamma$ are not mentioned among the conditions,
their terms which do not lead to another previously mentioned case 
can be cancelled in all of these cases.\\
To prove simplicity is rather easy for cases (d)-(f), because there are
no non-simple singularities of sufficiently small Tjurina number. For
cases (a)-(c), we obviously have adjacency relations (c) adjacent to
(b) with the same $k$, (b) to (a) again with the same $k$ and (a)
adjacent to (c) with a drop in $k$ by one. By lengthy, but explicit
calculations (using the 'adjacency'-relations among the 1-jets) one
can then rule out any other adjacencies than to the series of 
1-jet-types \JBE-\JBB and \JAA. As the latter ones are all simple,
this implies simplicity.
\item[\JAD-\JAF:] By the same kind of argument as for \JBA, no simple
singularities are possible.
\end{entry}
\end{proof}

\begin{rem} \label{tjulist}
As one can see in the last column of the table, we found precisely the
rational triple point singularities classified by Tjurina (see \cite{Tju}),
whose notation for their types is used. As the last three cases do
not have a name in the article of Tjurina, we simply stated them in the
same order as they appear there.
\end{rem}

\subsection{Singularities in $({\mathbb C}^5,0)$}

In this case we only need to consider matrices whose 1-jet involves at
least 4 variables. The methods are basically the same as in the
previous case, with one exception: For the case \JCA, the
problem of classification and of finding adjacencies can be reduced to
the corresponding problem for plane curve singularities and
deformations with sections thereof.

\begin{thm}
The simple isolated Cohen-Macaulay codimension 2 singularities in
$(\C^5,0)$ are the following ones:
\begin{center}
\begin{longtable}{|c|c|l|c|c|l|}
\hline
    Jet-Type
    & Type
    & Presentation Matrix
    &
    & $\tau$
\\ \hline
\JCB
    & $A^\sharp_0$
    & $\begin{pmatrix}x & y & z \\ w & v & x \end{pmatrix}$
    &
    & 1
\\ \hline
\JCA
    & $A^\sharp_k$
    & $\begin{pmatrix}x & y & z \\ w & v & x^{k+1}+y^2 \end{pmatrix}$
    & $k\geq 1$
    & $k+2$
\\ \hline
    & $D^\sharp_k$
    & $\begin{pmatrix}x & y & z \\ w & v & xy^2+x^{k-1} \end{pmatrix}$
    & $k\geq 4$
    & $k+2$
\\ \hline
    & $E^\sharp_6$
    & $\begin{pmatrix}x & y & z \\ w & v & x^3+y^4 \end{pmatrix}$
    &
    & 8
\\ \hline
    & $E^\sharp_7$
    & $\begin{pmatrix}x & y & z \\ w & v & x^3+xy^3 \end{pmatrix}$
    &
    & 9
\\ \hline
    & $E^\sharp_8$
    & $\begin{pmatrix}x & y & z \\ w & v & x^3+y^5 \end{pmatrix}$
    &
    & 10
\\ \hline
\JBE
    & $\Pi_k$
    & $\begin{pmatrix}w & y & x \\ z & w & y+v^k  \end{pmatrix}$
    & $k\geq 2$
    & $2k-1$
\\ \hline
\JBD
    &
    & $\begin{pmatrix} w & y & x  \\ z & w & y^k+v^2  \end{pmatrix}$
    & $k\geq 2$
    & $k+2$
\\ \hline
    &
    & $\begin{pmatrix} w & y & x \\ z & w & yv+v^k \end{pmatrix}$
    &
    & $2k$
\\
    &
    & $\begin{pmatrix} w+v^k & y & x \\ z & w & yv \end{pmatrix}$
    &
    & $2k+1$
\\ \hline
    &
    & $\begin{pmatrix} w+v^2 & y & x \\ z & w & y^2+v^k \end{pmatrix}$
    &
    & $k+3$
\\ \hline
    &
    & $\begin{pmatrix} w & y & x \\ z & w & y^2+v^3 \end{pmatrix}$
    &
    & $7$
\\ \hline
\JBC
    &
    & $\begin{pmatrix}v^2+w^k & y & x \\ z & w & v^2+y^l \end{pmatrix}$
    & $l\geq k\geq 2$
    & $k+l+1$
\\ \hline
    &
    & $\begin{pmatrix}v^2+w^k & y & x \\ z & w & yv \end{pmatrix}$
    & $k\geq 2$
    & $k+4$
\\ \hline
    &
    & $\begin{pmatrix}v^2+w^k & y & x \\ z & w & y^2+v^l \end{pmatrix}$
    & $k\geq 2, l\geq 3$
    & $k+l+2$
\\ \hline
    &
    & $\begin{pmatrix}wv+v^k & y & x \\ z & w & yv+v^k \end{pmatrix}$
    & $k\geq 3$
    & $2k+1$
\\
    &
    & $\begin{pmatrix}wv+v^k & y & x \\ z & w & yv \end{pmatrix}$
    & $k\geq 3$
    & $2k+2$
\\ \hline
    &
    & $\begin{pmatrix}wv+v^3 & y & x \\ z & w & y^2+v^3 \end{pmatrix}$
    &
    & $8$
\\ \hline
    &
    & $\begin{pmatrix}wv & y & x \\ z & w & y^2+v^3 \end{pmatrix}$
    &
    & $9$
\\ \hline
    &
    & $\begin{pmatrix}w^2+v^3 & y & x \\ z & w & y^2+v^3 \end{pmatrix}$
    &
    & $9$
\\ \hline
\JBB
    &
    & $\begin{pmatrix}z & y & x \\ x & w & v^2+y^2+z^k \end{pmatrix}$
    & $k\geq 2$
    & $k+4$
\\ \hline
    &
    & $\begin{pmatrix}z & y & x \\ x & w & v^2+yz+y^kw \end{pmatrix}$
    & $k\geq 1$
    & $2k+5$
\\
    &
    & $\begin{pmatrix}z & y & x \\ x & w & v^2+yz+y^{k+1} \end{pmatrix}$
    & $k\geq 2$
    & $2k+4$
\\ \hline
    &
    & $\begin{pmatrix}z & y & x \\ x & w & v^2+yw+z^2 \end{pmatrix}$
    &
    & 8
\\ \hline
    &
    & $\begin{pmatrix}z & y & x \\ x & w & v^2+y^3+z^2 \end{pmatrix}$
    &
    & 9
\\ \hline
    &
    & $\begin{pmatrix}z & y & x+v^2 \\ x & w & vy+z^2 \end{pmatrix}$
    &
    & 7
\\ \hline
    &
    & $\begin{pmatrix}z & y & x+v^2 \\ x & w & vz+y^2 \end{pmatrix}$
    &
    & 8
\\ \hline
    &
    & $\begin{pmatrix}z & y & x+v^2 \\ x & w & y^2+z^2 \end{pmatrix}$
    &
    & 9
\\ \hline
\end{longtable}
\end{center}
\end{thm}

\begin{proof} As in the case of 4 variables, we need to check all
possible 1-jets and exclude the non-simple singularities. To this end,
we start again by giving the list of those non-simple singularities
which we will need in the proof:

\begin{center}
\begin{longtable}{|c|l|c|c|c|}
\hline
    Jet-Type
    & Presentation Matrix
    & $\underline{a}$
    & $D$
    & $\tau$
\\ \hline
    \JCA
    & $\begin{pmatrix} x & y & z \\ w & v & x^4+y^4+\alpha \end{pmatrix}$
    & $\begin{pmatrix} 1 & 1 & 4 \\ 1 & 1 & 4 \end{pmatrix}$
    & $\begin{pmatrix} 1 & 1 & 4 & 1 & 1 \end{pmatrix}$
    & 11
\\ \hline
    & $\begin{pmatrix} x & y & z \\ w & v & x^3+y^6+\alpha \end{pmatrix}$
    & $\begin{pmatrix} 2 & 1 & 6 \\ 2 & 1 & 6 \end{pmatrix}$
    & $\begin{pmatrix} 2 & 1 & 6 & 2 & 1 \end{pmatrix}$
    & 12
\\ \hline
    \JBD
    & $\begin{pmatrix} w+v^2 & y & x \\ z & w & y^3+v^3 \end{pmatrix}$
    & $\begin{pmatrix} 2&1&2\\3&2&3 \end{pmatrix}$
    & $\begin{pmatrix} 2&1&3&2&1 \end{pmatrix}$
    & 8
\\ \hline
    & $\begin{pmatrix} w+v^3 & y & x \\ z & w & y^2+v^4 \end{pmatrix}$
    & $\begin{pmatrix} 3 & 2 & 3 \\ 4 & 3 & 4  \end{pmatrix}$
    & $\begin{pmatrix} 3 & 2 & 4 & 3 & 1 \end{pmatrix}$
    & 9
\\ \hline
    \JBB
    & $\begin{pmatrix}z & y & x \\ x & w & v^2+y^3+z^3(+yz^2+yw) \end{pmatrix}$
    &    $\begin{pmatrix}2 & 2 & 4 \\ 4 & 4 & 6\end{pmatrix}$
    &    $\begin{pmatrix}4 & 2 & 2 & 4 & 3 \end{pmatrix}$
    & 11
\\ \hline
    &    $\begin{pmatrix}z & y & x \\ x & w & v^3+y^2+z^2 \end{pmatrix}$
    &    $\begin{pmatrix} 2 & 3 & 4 \\ 4 & 5 & 6 \end{pmatrix}$
    &    $\begin{pmatrix} 4 & 3 & 2 & 5 & 2 \end{pmatrix}$
    & 13
\\ \hline
    &    $\begin{pmatrix}z & y & x \\ x & w & v^3+y^3+z^2 \end{pmatrix}$
    &    $\begin{pmatrix} 6 & 4 & 9 \\ 9 & 7 & 12 \end{pmatrix}$
    &    $\begin{pmatrix} 9 & 4 & 6 & 7 & 4 \end{pmatrix}$
    & 17
\\ \hline
    &    $\begin{pmatrix}z & y & x \\ x & w & v^2+y^4+z^2 \end{pmatrix}$
    &    $\begin{pmatrix}2 & 1 & 3 \\ 3 & 2 & 4 \end{pmatrix}$
    &    $\begin{pmatrix}3 & 1 & 2 & 2 & 2\end{pmatrix}$
    & 12
\\ \hline
    &    $\begin{pmatrix}z & y & x+v^2 \\ x & w & vz+yz+vw \end{pmatrix}$
    & $\begin{pmatrix} 3 & 2 & 4 \\ 4 & 3 & 5 \end{pmatrix}$
    & $\begin{pmatrix} 4 & 2 & 3 & 3 & 2 \end{pmatrix}$
    & 10
\\ \hline
    &    $\begin{pmatrix}z & y & x+v^3 \\ x & w & vy+z^2 \end{pmatrix}$
    & $\begin{pmatrix} 2 & 3 & 3 \\ 3 & 4 & 4 \end{pmatrix}$
    & $\begin{pmatrix} 3 & 3 & 2 & 4 & 1 \end{pmatrix}$
    & 9
\\ \hline
    &    $\begin{pmatrix}z & y & x+v^3 \\ x & w & y^2+yz+z^2 \end{pmatrix}$
    & $\begin{pmatrix} 2 & 2 & 3 \\ 3 & 3 & 4 \end{pmatrix}$
    & $\begin{pmatrix} 3 & 2 & 2 & 3 & 1 \end{pmatrix}$
    & 15
\\ \hline
    &    $\begin{pmatrix}z & y & x+v^2 \\ x & w & vy+yz+z^3 \end{pmatrix}$
    & $\begin{pmatrix} 1 & 2 & 2 \\ 2 & 3 & 3 \end{pmatrix}$
    & $\begin{pmatrix} 2 & 2 & 1 & 3 & 1 \end{pmatrix}$
    & 8
\\ \hline
    \JBA
    &    $\begin{pmatrix}x & y & \alpha \\ z & w & \beta \end{pmatrix}$
    &    $\begin{pmatrix}1 & 1 & 2 \\ 1 & 1 & 2\end{pmatrix}$
    &    $\begin{pmatrix} 1 & 1 & 1 & 1 & 1 \end{pmatrix}$
    & 17
\\ \hline
    \JBF
     & $\begin{pmatrix} x & y & z \\ w & \alpha & \beta \end{pmatrix}$
     & $\begin{pmatrix}1 & 1 & 1 \\ 2 & 2 & 2 \end{pmatrix}$
     & $\begin{pmatrix} 1 & 1 & 1 & 2 & 1 \end{pmatrix}$
     & 13
\\ \hline
\end{longtable}
\end{center}

\begin{entry}
\item[\JCB:] The matrix $j_1M$ already defines an isolated
singularity and is 1-determined, that is we may w.l.o.g assume that
$M=j_1M$. $M$ can be deformed to  
$$M'=\begin{pmatrix}x+\varepsilon & y & v \\ 
                   z & w & x\end{pmatrix}$$ 
with $\varepsilon\in\C$. Since the ideal defining the singular locus
of $M'$ is $\left<\varepsilon^2,\varepsilon v,\varepsilon
w,\varepsilon z,\ldots\right>$, the singular locus is empty, $M'$ is
smooth and hence $M$ simple. Because $M$ is contact-equivalent to
$$\begin{pmatrix}x & y & v \\ 
                 z & w & x+y^2\end{pmatrix},$$ 
we will call it $A^+_0$; the reason for this will become clear in the
subsequent case.\\

\item[\JCA:] Any matrix of type \JCA is contact-equivalent to a matrix
$$M\sim_C\begin{pmatrix}x & y & v \\ z & w & f(x,y)\end{pmatrix}$$ with $f(x,y)\in\m^2$.
We will show that the properties and the behaviour of the singularity defined by $M$ are determined by the hypersurface singularity defined by $f(x,y)$ in
$\C\{x,y\}$.

\begin{entry}
\item[Singular locus: ]By direct computation, the singular locus of
$M$ is completely contained in the plane defined by $\left<z,v,w\right>$, and in this plane, the singular locus of $M$ contains exactly the same points as the
singular locus of the singularity defined by $f$.

\item[Contact-Equivalence: ]If $f\in\C\{x,y\}$ is contact-equivalent
to some $g\in\C\{x,y\}$, there is an isomorphism $\gamma$ of $\C\{x,y\}$ and an unit $u\in\C\{x,y\}$ such that $ug=f\circ \gamma$. If $\gamma$ is given by
$\gamma(x)=\alpha_1x+\beta_1y$ and $\gamma(y)=\alpha_2x+\beta_2y$ with $\alpha_1,\beta_1,\alpha_2,\beta_2\in\C\{x,y\}$, we can extend $\gamma$ to an isomorphism of
$\C\{x,y,z,w,v\}$ by defining $\gamma(v)=uv$, $\gamma(z)=\alpha_1z+\beta_1w$ and $\gamma(w)=\alpha_2z+\beta_2w$. This is an isomorphism showing $M=\begin{pmatrix}x
& y & v \\ z & w & f(x,y)\end{pmatrix}$ is contact-equivalent to $\begin{pmatrix}x & y & v \\ z & w & g(x,y)\end{pmatrix}$.

\item[ Adjacencies: ]By direct computation, $T^1(M)$ can contain
only elements of the form
$$\begin{pmatrix}0&0&0\\0&0&h(x,y)\end{pmatrix}$$
with $h\in\C\{x,y\}\left/\left(f,x\frac{\partial f}{\partial x},y\frac{\partial f}{\partial x},x\frac{\partial f}{\partial y},y\frac{\partial f}{\partial
y}\right)\right.=\C\{x,y\}/(f,\m J(f))$. This is just the $T^1$ with section of the hypersurface singularity defined by $f$, and hence the adjacencies of $M$ are
determined by the adjacencies of $f$.
\end{entry}
In this way we get the simple isolated singularites $A^+_k, D^+_k, E^+_6, E^+_7$ and $E^+_8$. All of them can be deformed into the singularity $A^+_0$ we will get
in case \JCB corresponding to the smooth curve $f(x,y)=x+y^2$, $A_0$.

\item[\JBE:] Any matrix $M$ of type \JBE is contact-equivalent to a matrix
$$\begin{pmatrix} x & y & z+\delta_1v^{k_1} \\ 
                  z & w & y+\delta_2v^{k_2} \end{pmatrix}$$ 
with $\delta_1,\delta_2\in\{0,1\}$ and $k_1,k_2\geq 2$. 
If $\delta_1=\delta_2=0$, the singularity defined by $M$ is not
isolated. If $\delta_1\neq 0$ or $\delta_2\neq 0$, $M$ is contact-equivalent
to a matrix 
$$\begin{pmatrix} x & y & z \\ 
                  z & w & y+v^k \end{pmatrix}$$
with $k=\min\{k_i|\delta_i\neq 0\}$.(If $k_1 = k_2$, one of the
two terms may be cancelled by a lengthy sequence of coordinate changes
and row and column operations.)

We will call the singularity defined by $M$ $\Pi_k$. The singularity
$\Pi_k$ can only be deformed to 
$$M'\sim_c
  \begin{pmatrix} x & y & z+\alpha_{k-1}v^{k-1}+\ldots+\alpha_1 v+\alpha_0 \\
  z & w & y+v^k+\beta_{k-2}v^{k-2}+\ldots+\beta_1
                    v+\beta_0\end{pmatrix}$$ 
with $\alpha_i,\beta_j\in\C$, which is contact-equivalent to 
$$\begin{pmatrix} x & y & z \\
                  z & w & y+v^{k'} \end{pmatrix}$$ 
with $k'<k$. For $k'=0$, $M'$ is smooth, for $k'=1$, it is
contact-equivalent to the singularity $A^+_0$, and for $k'>1$, it is 
just $\Pi_{k'}$.

\item[\JBD:] A matrix of this type is contact-equivalent to a matrix
of the structure
$$\begin{pmatrix}w+\alpha & y & x\\
                 z & w & \beta\end{pmatrix}$$
where $\alpha,\beta \in {\mathfrak m}^2$, $\alpha$ only involving
terms in $v$, $\beta$ only involving terms in $y$ and $v$. By the
table of non-simple singularities (\JBD, line 1), simple singularities
cannot occur if the order of $\beta$ is at least 3. Therefore we may
assume (after a coordinate change in $v$ and cleaning up $w$ by column
operations) that the 2-jet of $\beta$ is one of the following:
$y^2+v^2$, $v^2$, $yv$ and $y^2$. In the first case, all terms in
$\alpha$ and the terms of higher order in $\beta$ may be cancelled
due to determinacy, and we obtain the first matrix of the first
series. If the 2-jet of $\beta$ is $v^2$, we can get rid of all terms
of $\alpha$ by substracting appropriate multiples of the last column
from the second one and then cleaning up by a sequence of coordinate
changes in $y$, $w$ and $z$ and a row operation. Since an appropriate
coordinate change in $v$ cancels all higher order terms in $\beta$
which involve $v$, the matrix in this case is (by determinacy) of the
structure 
$$\begin{pmatrix}w & y & x\\
                 z & w & v^2+y^k\end{pmatrix}$$
for some $k>2$, which is the first series. If the term $v^2$ is not
present, but the term $yv$ occurs, then all terms of higher order
which involve $y$ may be cancelled by a coordinate change in $v$ (and
possibly cleaning up again). This provides us with a matrix 
$$\begin{pmatrix}w+a_1v^{k_1} & y & x\\
                 z & w & yv+a_2v^{k_2}\end{pmatrix}$$
where $a_1,a_2 \in \C\{v\}$ either zero or a unit. By determinacy we
then obtain the second series. In the last of the four cases, we see
from the table of non-simple singularities (\JBD, line 2) that the
singularity cannot be simple if $\alpha$ is of order at least 3 and
the 3-jet of $\beta$ does not involve $v^3$. In other words, we can
only have further candidates for simple singularities, if the matrix
is of the structure
$$\begin{pmatrix}w+v^2 & y & x\\
                 z & w & y^2+\gamma\end{pmatrix}$$
or
$$\begin{pmatrix}w+\alpha & y & x\\
                 z & w & y^2+v^3\end{pmatrix}$$
where $\alpha,\gamma \in {\mathfrak m}^3$, $\alpha$ only involving $v$
and $\gamma$ only involving terms $v^l$ and $yv^l$. The terms $yv^l$
of $\gamma$ may be cancelled by appropriate coordinate changes and row
and column operations, since $l \geq 2$; this gives rise to the third
series. Again by sequence of column operations and coordinate changes
in $w$ and $z$, we can get rid of $\alpha$, providing the last
singularity for this case.\\
To proof simplicity, we can apply the same reasoning as in the case of
4 variables and obtain that, in addition to the adjacencies which
preserve the 1-jet of the matrix, adjacencies are at most possible to
singularities of the $A$-series, of 1-jet \JBE, \JCB and
\JCA. For adjacencies preserving the 1-jet of the matrix, we may again
follow the classification to see that adjacencies to the non-simple
singularities of this 1-jet are impossible. A direct but rather
lengthy calculation (relying on the fact that all entries are of order
at most 2 in the original matrix and tracing this fact throughout the
whole computation) then shows that no adjacencies to non-simple
singularities of type \JCA are possible. 

\item[\JBC:] In this case, the matrix is of the structure \footnote{At
this point it is important to observe that the roles of $y$ and $w$
may harmlessly be interchanged.}  
$$\begin{pmatrix} \alpha & y & x\\
                  z & w & \beta\end{pmatrix}$$
where $\alpha,\beta \in {\mathfrak m}^2$, $\alpha$ only involving $w$
and $v$, $\beta$ only involving $y$ and $v$. By the table of
non-simple singularities (\JBD, line 2), we see immediately
that no simple singularities can occur if the order of $\alpha$ or
$\beta$ exceeds 2. As $\alpha$ and $\beta$ both have to be of order 2,
it turns out to be a suitable approach to consider three cases: 
\begin{enumerate}
\item[(a)] $\alpha$ and $\beta$ contain a term $v^2$\\
By appropriate coordinate changes in $v$ followed by applying column
operations of type 'addition of monomial times the second column to the
first' resp. 'to the third column' and suitable 
coordinate changes in $x$ and $z$, we can obtain a matrix of the
following structure after applying the usual determinacy argument:
$$\begin{pmatrix}v^2+w^k & y & x\\
                 z & w & v^2+y^l\end{pmatrix}$$
These are exactly the matrices of the first series.
\item[(b)] only one of $\alpha$ and $\beta$ contains a term $v^2$\\
W.l.o.g. we may assume that the term $v^2$ appears in $\alpha$. By the
same transformation as in the case (a), the matrix is then of the
structure
$$\begin{pmatrix}v^2+w^k & y & x\\
                 z & w & \beta\end{pmatrix}$$
where $\beta$ only involves $v$ and $y$, is of order 2 and does not
contain the term $v^2$. Therefore the 2-jet of $\beta$ is of the form
$ay^2+byv$. If $b\neq 0$, we can kill the term $ay^2$ as well
as all higher pure powers of $y$ by a suitable coordinate change in
$v$ (and of course subsequent cleaning up in the upper left hand
entry). By applying a determinacy argument the matrix is then of the
structure
$$\begin{pmatrix}v^2+w^k & y & x\\
                 z & w & yv\end{pmatrix}.$$
If $b=0$, $\beta$ is of the structure $\tilde{a}y^2+cyv^2+dv^l$ where
$\tilde{a},b,c \in \C\{y,v\}$, $\tilde{a}$ and $d$ units. We can get
rid of the $cyv^2$ by an appropriate coordinate change in $y$,
applying column operations of type 'addition of monomial times first
column to the second' and subsequent cleanup by row and column
operations and coordinate changes in $x$ and $z$. After applying a
determinacy argument, this provides us with a matrix of the structure
$$\begin{pmatrix}v^2+w^k & y & x\\
                 z & w & y^2+v^l\end{pmatrix}$$
where $l\geq 3$ which is exactly the third series.
\item[(c)] neither $\alpha$ nor $\beta$ contain a term $v^2$ \\
By the same arguments as at the beginning of case (b), we need to
distinguish between three cases: the 2-jets of $\alpha$ and $\beta$
both still involve $v$, only one involves $v$ or both do not involve
$v$. In the first case, the 2-jet of the matrix is
$$\begin{pmatrix}wv & y & x\\
                 z & w & yv\end{pmatrix}$$
and the only higher order terms in $\alpha$ and $\beta$ which cannot
be cancelled are pure powers of $v$ yielding exactly the 4th series.\\
If there is only one mixed term in the 2-jet, we may w.l.o.g. assume
that it appears in $\alpha$. By the same kind of sequence of row and
column operations and coordinate changes as before we obtain a matrix
of the structure
$$\begin{pmatrix}wv+av^3 & y & x\\ 
                 z & w & y^2+bv^3+cyv^2\end{pmatrix}$$
where $a,b,c \in \C\{v\}$. If $b$ is not a unit, perturbing the upper
left hand entry with $w$ leads to the non-simple singularity (\JBD,
line 2). By an appropriate coordinate change in $v$ (and subsequent
cleaning up), we may hence safely assume that $b=1$ and $c=0$. If $a$
is a unit, we are dealing with the first of the three remaining
singularities.\\
If $a$ is not a unit, we can write $av^3$ as
$\tilde{a}v^4$ and a coordinate change in $w$ can move it to the
middle entry on the bottom row where it appears as $\tilde{a}v^3$
which can in turn be cancelled by adding $\tilde{a}$ times the third
column to the second one (and subsequent cleaning up). This
singularity is the second one of the remaining three.\\
In the final case, where $\alpha$ and $\beta$ do not involve $v$, the
only singularity which is not adjacent to the non-simple one already
mentioned above is
$$\begin{pmatrix}w^2+v^3 & y & x\\
                 z & w & y^2+v^3\end{pmatrix}$$
which is exactly the last singularity in the list.
\end{enumerate}
The proof that these singularities are indeed simple involves the same
kind of arguments as in the case \JBD, allowing, in addition to
adjacencies preserving the 1-jet, at most adjacencies to singularities
of the $A$-series, of 1-jet \JCB, \JCA, \JBD and \JBE. By basically
the same lengthy calculations as in the case \JBD, it can again be
shown that adjacencies to the non-simple singularities of the same
1-jet, of 1-jet \JBD and of 1-jet \JCA cannot occur.

\item[\JBB:] A matrix with 1-jet of type \JBB is of the structure
$$\begin{pmatrix}z & y & x+\alpha\\
                 x & w & \beta\end{pmatrix}$$
where $\alpha,\beta \in {\mathfrak m}^2$, $\alpha$ only involving
$v$. From the table of non-simple singularities (\JBB, line1),
we see first of all that simple singularities can only occur if the
order of $\beta$ is 2 and by (\JBB, all lines) that for simple
singularities the 2-jet of $\beta$ cannot consist of a single pure
power. Moreover, direct computation of the singular locus shows that
we need at least one pure power of $v$ in $\alpha$ or $\beta$.
This still leaves a rather large number of possibilities for
the 2-jet of $\beta$ which we can divide into four cases: 
\begin{enumerate}
\item[(a)] $\beta$ contains the term $v^2$\\
By a suitable coordinate change in $v$, we may assume that $v^2$ is
the only term in $\beta$ which involves $v$. Using the fact that the
2-jet of $\beta$ cannot consist of a single pure power, a direct 
calculation shows that the following 7 cases may occur (which we
already saw in the case of 4 variables): $v^2+y^2+z^2$, $v^2+yz+yw$,
$v^2+z^2+yw$, $v^2+y^2$, $v^2+yz$, $v^2+z^2$, $v^2+yw$. In the first
three cases the matrix is already 2-determined and due to weighted
determinacy we can even get rid of all terms of $\alpha$ obtaining
three singularities of which the first two are just the beginning of
the first two series. The remaining cases give rise to two more series
and two additional singularities by the same reasoning as in the case
\JBB in 4 variables each time, of course, using the weighted
determinacy to remove the terms of $\alpha$.
\item[(b)] $\beta$ does not contain $v^2$, but $vy$\\
In this case, direct calculation shows that we can have the following
2-jets of $\beta$: $vy+yz+z^2$,  $vy+yz$, $vy+z^2$ and $vy$. If the
order of $\alpha$ is at least three we can perturb a matrix in this
case to the first of the 4 possibilities for $\beta$ and to $v^3$ for
$\alpha$; in this matrix the term $yz$ of $\beta$ is then killed by
determinacy and we obtain the matrix \JBB, line 6, from the table of
non-simple singularities. Hence simple singularities can only occur if
$\alpha$ contains the term $v^2$. For determinacy reasons, the first
and the third of the above 4 cases then coincide leading to the
singularity
$$\begin{pmatrix}z & y & x+v^2\\
                 x & w & vy+z^2\end{pmatrix}.$$
The other two cases cannot lead to simple singularities according to
the list of non-simple singularities.
\item[(c)] $\beta$ contains neither $v^2$ nor $vy$, but $vz$\\
Similar to the previous case, the only possible 2-jets of $\beta$ are
$vz+y^2+yz$, $vz+y^2$, $vz+yz$ and $vz$. As perturbing $\beta$ with
the term $vy$ takes us to the case (b), we see that $\alpha$ can be
of order at most 2. By determinacy, the first and second case then
give rise to the same singularity
$$\begin{pmatrix}z & y & x+v^2\\
                 x & w & vz+y^2\end{pmatrix}.$$
The other two cases are adjacent (by the perturbation term $vy$) to
the non-simple singularity mentioned in (b).
\item[(d)] $\beta$ does not contain any terms involving $v$\\
For the 2-jet of $\beta$, a direct calculation with row and column
operations and coordinate changes in $y$ and $w$ shows that there are
not many possibilities left in this case: $y^2+z^2$, $y^2$, $yz$ and
$z^2$. In all of these cases, the singularities are adjacent to
singularities from case (b) which implies that $\alpha$ has to be of
order 2 for simple singularities. Moreover, the second case is
adjacent to the non-simple singularity from case (b) (by perturbation
with $vy$) and the last two are adjacent to non-simple singularities
from case (c) (by perturbation with $vz$). This only leaves the first
case, namely
$$\begin{pmatrix}z & y & x+v^2\\
                 x & w & y^2+z^2\end{pmatrix},$$
which is the last singularity in the list stated in the theorem. 
\end{enumerate}
For proving simplicity of these singularities, it is first of all
important to show that the singularities, which we found as
candidates in this case \JBB, cannot be adjacent to any of the
non-simple singularities of the same case. Due to the large number of
cases, this is a lengthy, but straightforward calculation. On the
other hand, adjacencies to singularities with other 1-jets of the
matrix can be checked explicitly, which leads possible adjacencies to  
the D-series, the A-series, singularities of 1-jet \JBE, \JBD, \JBC,
\JCB and \JCA. The main ingredients to computing directly that
adjacencies to non-simple singularities of types \JCA, \JBC and \JBD
cannot occur, are the knowledge about the low order of $\alpha$ and
$\beta$ and about the structure of $\beta$ in each of the cases;
i.e. the 4 cases (a) to (d) are considered separately which makes the
calculation even lengthier then the previous ones.

\item[\JBA:] According to the table of non-simple singularities the
generic matrix of this 1-jet cannot be simple.

\item[\JBF:] Again there cannot be any simple singularities according
to the table of non-simple singularities.
\end{entry}
\end{proof}

\subsection{Singularities in $({\mathbb C}^6,0)$}
In this case, we can only exclude the 1-jets containing 3 variables
which implies that we have to consider jets containing 4, 5 and 6
variables. Fortunately, there turns out to be only one case containing
6 variables and only very few cases with 4 variables contributing to the
list of simple singularities. Parts of the proof parallel the
classifications of simple hypersurface singularities of dimension 2,
other parts rely on the one of simple fat point singularities in the
plane. 

\begin{thm}
The simple isolated Cohen-Macaulay codimension 2 singularities in
$(\C^6,0)$ are listed in the following table:.
\begin{center}
\begin{longtable}{|c|c|l|c|c|l|}
\hline
    Jet-Type
    & Type
    & Presentation Matrix
    &
    & $\tau$
\\ \hline
\JDA
    & $\Omega_1$
    & $\begin{pmatrix}x & y & v \\ z & w & u \end{pmatrix}$
    &
    & 0
\\ \hline
\JCB
    & $\Omega_k$
    & $\begin{pmatrix}x & y & v \\ z & w & x+u^k \end{pmatrix}$
    & $k\geq 2$
    & $k-1$
\\ \hline
\JCA
    & $A_k^\sharp$
    & $\begin{pmatrix}x & y & z \\ w & v & u^2+x^{k+1}+y^2 \end{pmatrix}$
    & $k\geq 1$
    & $k+2$
\\ \hline
    & $D_k^\sharp$
    & $\begin{pmatrix}x & y & z \\ w & v & u^2+xy^2+x^{k-1} \end{pmatrix}$
    & $k\geq 4$
    & $k+2$
\\ \hline
    & $E_6^\sharp$
    & $\begin{pmatrix}x & y & z \\ w & v & u^2+x^3+y^4 \end{pmatrix}$
    &
    & 8
\\ \hline
    & $E_7^\sharp$
    & $\begin{pmatrix}x & y & z \\ w & v & u^2+x^3+xy^3 \end{pmatrix}$
    &
    & 9
\\ \hline
    & $E_8^\sharp$
    & $\begin{pmatrix}x & y & z \\ w & v & u^2+x^3+y^5 \end{pmatrix}$
    &
    & 10
\\ \hline
    &
    & $\begin{pmatrix}x & y & z \\ w & v & ux+y^k+u^l \end{pmatrix}$
    & $k\geq 2, l\geq 3$
    & $k+l-1$
\\ \hline
    &
    & $\begin{pmatrix}x & y & z \\ w & v & x^2+y^2+u^3 \end{pmatrix}$
    &
    & 6
\\ \hline
\JBE
    & $F_{q,r}^\sharp$
    & $\begin{pmatrix}w & y & x \\ z & w+vu & y+v^q+u^r \end{pmatrix}$
    & $q,r\geq 2$
    & $q+r$
\\ \hline
    & $G_5^\sharp$
    & $\begin{pmatrix}w & y & x \\ z & w+v^2 & y+u^3 \end{pmatrix}$
    &
    & 7
\\ \hline
    & $G_7^\sharp$
    & $\begin{pmatrix}w & y & x \\ z & w+v^2 & y+u^4 \end{pmatrix}$
    &
    & 10
\\ \hline
    & $H_{q+3}^\sharp$
    & $\begin{pmatrix}w & y & x \\ z & w+v^2+u^q & y+vu^2 \end{pmatrix}$
    & $q\geq 3$
    & $q+5$
\\ \hline
    & $I_{2q-1}^\sharp$
    & $\begin{pmatrix}w & y & x \\ z & w+v^2+u^3 & y+u^q \end{pmatrix}$
    & $q\geq 4$
    & $2q+1$
\\ \hline
    & $I_{2r+2}^\sharp$
    & $\begin{pmatrix}w & y & x \\ z & w+v^2+u^3 & y+vu^r \end{pmatrix}$
    & $r\geq 3$
    & $2r+4$
\\ \hline
\JBD
    &
    & $\begin{pmatrix}w & y & x \\ z & w+v^{k_1}+u^{k_2} & y^l+uv \end{pmatrix}$
    & $k_1,k_2,l\geq 2$
    & $k_1+k_2+l-1$
\\ \hline
    &
    & $\begin{pmatrix}w & y & x \\ z & w+v^2 & u^2+yv \end{pmatrix}$
    &
    & 6
\\ \hline
    &
    & $\begin{pmatrix}w & y & x \\ z & w+uv & u^2+yv+v^k \end{pmatrix}$
    & $k\geq 3$
    & $k+4$
\\ \hline
    &
    & $\begin{pmatrix}w & y & x \\ z & w+v^k & u^2+yv+v^3 \end{pmatrix}$
    & $k\geq 3$
    & $2k+2$
\\
    &
    & $\begin{pmatrix}w & y & x \\ z & w+uv^k & u^2+yv+v^3 \end{pmatrix}$
    & $k\geq 2$
    & $2k+5$
\\ \hline
    &
    & $\begin{pmatrix}w & y & x \\ z & w+v^3 & u^2+yv \end{pmatrix}$
    &
    & 9
\\ \hline
    &
    & $\begin{pmatrix}w & y & x \\ z & w+v^k & u^2+y^2+v^3 \end{pmatrix}$
    & $k\geq 3$
    & $2k+3$
\\
    &
    & $\begin{pmatrix}w & y & x \\ z & w+uv^k & u^2+y^2+v^3 \end{pmatrix}$
    & $k\geq 2$
    & $2k+6$
\\ \hline
\end{longtable}
\end{center}
\end{thm}

\begin{proof}
First of all, we state a table of non-simple singularities in 6
variables:\\
\begin{center}
\begin{longtable}{|c|l|c|c|c|}
\hline
    Jet-Type
    & Presentation Matrix
    & $\underline{a}$
    & $D$
    & $\tau$
\\ \hline
    \JCA
    & $\begin{pmatrix} x & y & z \\ w & v & u^2+x^3+y^6 \end{pmatrix}$
    & $\begin{pmatrix} 2 & 1 & 6 \\ 2 & 1 & 6 \end{pmatrix}$
    & $\begin{pmatrix} 2 & 1 & 6 & 2 & 1 & 3\end{pmatrix}$
    & 12
\\ \hline
    & $\begin{pmatrix} x & y & z \\ w & v & u^2+x^4+y^4 \end{pmatrix}$
    & $\begin{pmatrix} 1 & 1 & 4 \\ 1 & 1 & 4 \end{pmatrix}$
    & $\begin{pmatrix} 1 & 1 & 4 & 1 & 1 & 2 \end{pmatrix}$
    & 11
\\ \hline
    & $\begin{pmatrix} x & y & z \\ w & v & x^2+y^2+\alpha \end{pmatrix}$
    & $\begin{pmatrix} 2 & 2 & 5 \\ 1 & 1 & 4 \end{pmatrix}$
    & $\begin{pmatrix} 2 & 2 & 5 & 1 & 1 & 1 \end{pmatrix}$
    & 8
\\ \hline
    & $\begin{pmatrix} x & y & z \\ w & v & y^2+x^3+u^3 \end{pmatrix}$
    & $\begin{pmatrix} 2 & 3 & 7 \\ 1 & 2 & 6 \end{pmatrix}$
    & $\begin{pmatrix} 2 & 3 & 7 & 1 & 2 & 2 \end{pmatrix}$
    & 8
\\ \hline
\JBE
    & $\begin{pmatrix} w & y & x \\ z & w+u^2 & y+v^5   \end{pmatrix}$
    & $\begin{pmatrix} 4 & 5 & 6 \\ 3 & 4 & 5 \end{pmatrix}$
    & $\begin{pmatrix} 6 & 5 & 3 & 4 & 1  & 2 \end{pmatrix}$
    & 13
\\ \hline
    & $\begin{pmatrix} w & y & x \\ z & w+u^3 & y+v^3   \end{pmatrix}$
    & $\begin{pmatrix} 3 & 3 & 3 \\ 3 & 3 & 3 \end{pmatrix}$
    & $\begin{pmatrix} 3 & 3 & 3 & 3 & 1 & 1 \end{pmatrix}$
    & 12
\\ \hline
\JBD
    & $\begin{pmatrix} w & y & x \\ z & w+\alpha & \beta   \end{pmatrix}$
    & $\begin{pmatrix} 4 & 3 & 3 \\ 5 & 4 & 4 \end{pmatrix}$
    & $\begin{pmatrix} 3 & 3 & 5 & 4 & 1 & 2 \end{pmatrix}$
    & 11
\\ \hline
    & $\begin{pmatrix} w & y & x \\ z & w+\alpha & \beta   \end{pmatrix}$
    & $\begin{pmatrix} 2 &2 & 3 \\2 & 2 & 3 \end{pmatrix}$
    & $\begin{pmatrix} 3 & 2 & 2 & 2 & 1 & 1 \end{pmatrix}$
    & 10
\\ \hline
    & $\begin{pmatrix} w & y & x \\ z & w+\alpha & \beta   \end{pmatrix}$
    & $\begin{pmatrix} 6 & 2 & 4 \\ 6 & 4 & 6 \end{pmatrix}$
    & $\begin{pmatrix} 4 & 2 & 6 & 4 & 2 & 3 \end{pmatrix}$
    & 9
\\ \hline
\end{longtable}
\end{center}

\begin{entry}
\item[\JDA:] By determinacy a singularity with this 1-jet is of the form
$$\begin{pmatrix} x & y & v \\ z & w & u\end{pmatrix}.$$
Since $T^1(M)=0$, this singularity is rigid and in particular simple.

\item[\JCB:] In this case, the matrix is of the structure
$$\begin{pmatrix} x & y & v\\ 
                 z & w & x+u^j\end{pmatrix}$$ 
for some  $j\geq 2$. Calculating $T^1(M)$, we see that the versal
deformation of $M$ is of the form 
$$\begin{pmatrix}x & y & v \\
               z & w &x+u^j+\sum\limits_{i=0}^{j-2}\alpha_iu^i\end{pmatrix}$$
where $\alpha_i\in\C$ and not all $\alpha_i=0$.
At each zero of $u^j+\sum\limits_{i=0}^{j-2}\alpha_iu^i$, this matrix
is contact-equivalent to the rigid singularity $\Omega_1$ or to a germ
$\Omega_k$ for some $k<l$. Thus all $\Omega_j$ for $j\geq 2$ are simple.

\item[\JCA:] Here the matrix is of the structure
$$\begin{pmatrix}x & y & z\\
                 w & v & \alpha\end{pmatrix}$$
where $\alpha \in {\mathfrak m}^2$, not involving $z$, $w$ and $v$.
In contrast to \JCA in 5 variables, we can not directly reduce this
case to the classification of simple isolated hypersurface
singularities in 3 variables, because the variable $u$ which does not
appear in the 1-jet of the matrix plays a different role than the
variables $x$ and $y$. Thus we will distinguish between three
cases\footnote{Note that the roles of $x$ and $y$ can harmlessly be
exchanged which explains why we can assume in the second case that the
mixed term is $ux$.} depending on the way $u$ appears in the 2-jet of
$\alpha$: 
\begin{enumerate}
\item[(a)] the 2-jet of $\alpha$ contains $u^2$\\
In this case, the other terms containing $u$ may be cancelled directly
by an appropriate coordinate change in $u$ leading to a matrix 
$$\begin{pmatrix}x & y & z\\
                 w & v & u^2+\beta\end{pmatrix}$$
where $\beta \in {\mathfrak m}^2$, only involving $x$ and $y$. By the
same calculations as in the case \JCA in 5 variables, it can now be
shown that the singular locus corresponds to
$(z,v,w,u,\beta,\frac{\partial \beta}{\partial x},
      \frac{\partial \beta}{\partial y})$ implying that points of the
singular locus correspond exactly to the points of the one of $\beta$.
The reformulation of contact equivalence and $T^1$ also lead to the
same observations as in the case of 5 variables. Therefore we obtain
simple singularities from this case exactly for $\beta$ being an $E_6$,
$E_7$, $E_8$, $D_k$ or $A_k$ singularity.

\item[(b)] the 2-jet of $\alpha$ contains $ux$, but not $u^2$\\
Here we can apply an appropriate coordinate change in $x$ (followed by
cleaning the first column by a suitable column operation of type
'adding the second column to the first' and a subsequent coordinate
change in $w$) and hence assume that there are no terms in $\alpha$
which are divisible by $yu$. A further coordinate change, this time in
$u$ allows us to get rid of all terms containing the factor $x$ except
$xu$, of course, yielding a matrix of the structure
$$\begin{pmatrix}x & y & z\\
                 w & v & ux+ay^k+bu^l\end{pmatrix}$$
where $a \in \C\{y\}$, $b \in \C\{u\}$ are both units, because a
matrix of this structure does not define an isolated singularity if
there is no pure power of $u$ or no pure power of $y$ in the lower
right hand entry. By determinacy, this is exactly the last series
listed for this 1-jet in the table of the theorem. 

\item[(c)] the 2-jet of $\alpha$ does not contain any terms involving $u$\\
In this case, we see from the table of non-simple singularities
(\JCA, line 4) that a simple singularity may only occur if the 2-jet of
$\alpha$ is not a square. This implies that we only need to consider
matrices of the structure
$$\begin{pmatrix}x & y & z\\
                 w & v & x^2+y^2+\beta\end{pmatrix}$$
where $\beta \in {\mathfrak m}^3$. According to the table of
non-simple singularities (\JCA, line 3), all singularities in this case
cannot be simple unless they contain the term $u^3$. In this case
determinacy implies that this is exactly the last singularity in the
list for this 1-jet.
\end{enumerate}

Since adjacencies to singularities of hypersurface types cannot occur in
this case and since all singularities of 1-jet \JDA and \JCA are simple,
we only need to consider adjacencies which do not change the 1-jet. In
the case (a) all of these adjacency calculations are exactly analogous
to the case of 2-dimensional hypersurfaces proving simplicity of the
respective singularities. In case (b) the only adjacencies preserving
the 1-jet are the ones into singularities of the same case and into
$A_k$ type of singularities of case (a), as a direct calculation
shows; for the last singularity we can already deduce from the Tjurina
number that adjacencies to non-simple singularities of the cases (a)
and (b) are impossible.

\item[\JBE:] A singularity of this kind corresponds to a matrix of the
structure
$$\begin{pmatrix}w & y & x\\
                 z & w+\alpha & y+\beta\end{pmatrix}$$
where $\alpha,\beta \in {\mathfrak m}^2$, involving only the variables
$u$ and $v$. Similar to the case \JCA in 5 variables, the best
approach to this case is to first consider the singular locus of these
singularities. It turns out to be defined by $(x,y,z,w,\alpha,\beta)$
implying that the given matrix can only define an isolated singularity
if $(\alpha,\beta)$ corresponds to a fat point in the plane. A
rather direct calculation shows that two singularities defined by
matrices of this kind are contact-equivalent if and only if the
corresponding fat points are contact-equivalent. Moreover, the $T^1$
can only contain elements of the form
$$\begin{pmatrix}0 & 0 & 0\\
                 0 & a & b\end{pmatrix}$$
where 
$$\begin{pmatrix}a\\b\end{pmatrix} \in (\C\{u,v\}/(\alpha,\beta))^2
 \left/
 \left(\begin{pmatrix}\frac{\partial \alpha}{\partial u}\\
                 \frac{\partial \beta}{\partial u}\end{pmatrix},
  \begin{pmatrix}\frac{\partial \alpha}{\partial v}\\
                 \frac{\partial \beta}{\partial v}\end{pmatrix}
 \right)\right..$$
This is exactly the $T^1$ of the corresponding fat point singularity.
Therefore simple singularities with the given 1-jet of the matrix are
exactly the ones of the above mentioned structure for which
$(\alpha,\beta)$ is a fat point from the list of simple isolated
complete intersection singularities which are not of type $A_k$. 

\item[\JBD:] For studying matrices in this case, the first important
observation is that these matrices are adjacent to those of 1-jet
\JBE by perturbing with $y$ in the lower right hand entry. Writing a
matrix of the given 1-jet as
$$\begin{pmatrix}w & y &x\\
                 z & w+\alpha & \beta\end{pmatrix}$$
we can hence deduce that this singularity can only be simple if at
least one of $\alpha$ and $\beta$ contains one of the order two terms
$uv$ and $u^2$ (or, of course, $v^2$ which coincides with the $u^2$
case by exchanging the roles of $u$ and $v$ which is still possible at
this point).
\begin{enumerate}
\item[(a)] $j_2\beta$ not a square modulo $y$\\
By appropriate row and column operations and coordinate changes of
$u$, $v$ and $z$, we may assume that a matrix in this case is of the
structure
$$\begin{pmatrix}w & y & x\\
                 z & w+a_1u^{k_1}+a_2v^{k_2} & uv+by^l\end{pmatrix}$$
for suitable $k_1,k_2,l$ and $a_1 \in \C\{u\}$, $a_2 \in \C\{v\}$
and $b \in \C\{y\}$ units, because the singularity is no longer
isolated if any of the three terms is missing. By determinacy, this
leads to the first series of \JBD.

\item[(b)] $j_2\beta=u^2$ modulo $y$\\
In this case, a coordinate change in $u$ allows us to cancel the term
$yu$ in $\beta$ and the term $u^2$ in $\alpha$ can be killed by a
sequence of appropriate row and column operations followed by a
coordinate change in $z$ (and subsequent cleanup). Therefore the 2-jet
of the matrix can be assumed to be of the form
$$\begin{pmatrix}w & y & x\\
                 z & w+auv+bv^2 & u^2+cyv+dy^2\end{pmatrix}.$$
where $a,b,c,d \in \C$. If $c \neq 0$, a coordinate change in $v$
followed by one in $w$, a column operation of type 'adding second
column to first column' and a coordinate change in $z$ allows us to
remove the term $dy^2$. If in addition to that also $b\neq 0$, then
determinacy yields that the matrix is
$$\begin{pmatrix}w & y & x\\
                 z & w+v^2 & u^2+yv\end{pmatrix}.$$
Still in the case $c\neq0$, but this time $b=0$, $a\neq0$, the
situation is a little bit more difficult, since we cannot use
determinacy to remove the higher order pure powers of $v$ in $\alpha$.
Instead, we make a coordinate change in $u$ which produces two kinds
of terms in $\beta$: terms of the form $uv^l$ which we will remove and
pure powers in $v$ which we will keep. For cancelling the terms $uv^l$,
we first perform a column operation of type 'adding second column to
the third' (and clean up the upper right hand entry by a coordinate
change in $x$) replacing $uv^l$ by $v^{l-1}w$. A coordinate change in
$y$ now allows us to move these terms to terms $v^{l-2}w$ in the upper
middle entry and by a column operation of type 'adding first column to
the second' to the bottom middle entry as $v^{l-2}z$. As $l$ was at
least three in our case, we can now shift these terms back to the
lower right hand entry by a suitable coordinate change in $u$ and
eventually kill them by a column operation of type 'adding first
column to third column'. By determinacy the matrix is then of the form
$$\begin{pmatrix}w & y & x\\
                 z & w+uv & u^2+yv+v^k\end{pmatrix}$$
for a suitable $k \geq 3$ (If the $v^k$ term were missing the
singularity would not stand a chance to be isolated.)\\
Again still in the case $c\neq0$, but now in the situation $a=b=0$, we
can conclude from the table of non-simple singularities (\JBD, all lines)
that a simple singularity can only occur if the term $v^3$ is present
in at least one of $\alpha$ and $\beta$. If it is present in $\beta$,
then the matrix is of the structure
$$\begin{pmatrix}w & y & x\\
                 z & w+\gamma_1 uv^{k_1} + \gamma_2 v^{k_2} & u^2+yv+v^3
  \end{pmatrix}$$
where $\gamma_1,\gamma_2 \in \C\{v\}$ and at least one of them a unit.
By determinacy, exactly one term $uv^{k_1}$ or $v^{k_2}$ remains and
we obtain the third series. If there is no $v^3$ in $\beta$,
then it appears in $\alpha$ leading, by determinacy, to the matrix
$$\begin{pmatrix}w & y & x\\
                 z & w+v^3 & u^2+yv\end{pmatrix}.$$
This ends the arguments in the case $c\neq 0$.\\
In the case $c=0$, $d \neq 0$, simple singularities may only occur if
the term $v^3$ is present in $\beta$, since otherwise an adjacency to
the non-simple singularity of jet type \JCA, line 3,
exists by perturbation of the middle entry in the bottom row by $u$.
If $v^3$ is present in $\beta$, then the matrix is of the structure
$$\begin{pmatrix}w & y & x\\
                 z & w+a_1uv^{k_1}+a_2v^{k_2} & u^2+y^2+v^3\end{pmatrix}$$
where $a_1,a_2 \in \C\{v\}$ either a unit or zero, at least one of
them a unit. In contrast to the previous cases, we cannot proceed by
the standard determinacy argument here, because some of these matrices
are not quasihomogeneous in the strict sense, as $x$ needs to be assigned
non-positive weights if $k_1$ exceeds 2 and $k_2$ exceeds $4$ in order
to satisfy relative row and column weights. On the other hand, we can
remove the higher order terms (w.r.t. the weights of $u^2+y^2+v^3$) by
explicit calculation in the following way: If the lowest order term is
of the form $uv^{k_1}$ we start by a coordinate change in $u$ killing
$a_2v^{k_2}$ which in turn introduces some $uv^2 \cdot p(v)$ into
the lower left hand entry of which all terms except at most one are
divisible by $v^3$ and can hence be collected into a coefficient of
$v^3$ which is a unit. We can get rid of the last remaining term by a
coordinate change in $v$ such that all new terms can again be
collected into the coefficients of $u^2$ and $v^3$. By appropriate
coordinate changes of $y$, $u$ and $w$ (multiplication by units), we can
achieve that the matrix is of the form 
$$\begin{pmatrix}e_1w & e_2y & x\\
                 z & e_1(w+uv^{k_1}) & e_3(u^2+y^2+v^3)\end{pmatrix}$$
where $e_1,e_2,e_3$ are units. These units can easily be removed by
multiplication of rows and columns by units and subsequent coordinate
changes in $z$ and $x$. If, on the other hand, the lowest order term
is of the form $v^{k_2}$, we start by a coordinate change in $v$ and
then remove the offending terms in the same way as before.\\
Finally, if $c=d=0$ in the matrix 
$$\begin{pmatrix}w & y & x\\
                 z & w+auv+bv^2 & u^2+cyv+dy^2\end{pmatrix}.$$
(from the beginning of this case) then the table of
non-simple singularities (\JBD, line 3) implies that there cannot be
any simple singularities. \\
For proving simplicity, the rather straightforward explicit
calculation shows in this situation that at most singularities of types
$A_k$, \JDA, \JCB, \JCA, \JBE and of course of the same 1-jet may
occur. $A_k$, \JDA and \JCB do not contain non-simple singularities; that
the non-simple ones of \JCA, \JBE and \JBD cannot be reached can be
shown directly by a calculation which relies on keeping track of the
order of the entries in $u$ and $v$.

\item[(c)] $j_2\beta=0$ modulo $y$\\
In this case, the table of non-simple singularities (\JBD, line 2)
implies that no simple singularities can occur.
\end{enumerate}

\item[\JBC-\JBF:] According to the table of non-simple singularities,
matrices with these 1-jets cannot define simple singularities.
\end{entry}
\end{proof}

\subsection{Some Remarks on the Lists of Simple Singularities}

In contrast to the table of simple singularities in 4 variables which
does not show any unexpected behaviour, there are several surprising
details about the tables in 5 and 6 variables:
\begin{itemize}
\item The simple hypersurface singularities reappear as entries of the
      matrix in the cases \JCA. This can easily be understood as a
      consequence of the fact that in the other 5 entries no
      perturbation exist which cannot be shifted into the last
      entry. Moreover only those variables cannot be removed from this
      last entry by row and column operations which are neither an
      entry in the same row nor in the same column.\\
\item The simple fat point singularities in the plane reappear as
      entries of the matrix in the case \JBE in 6 variables resp. the
      fat points on the line in the same case in 5 variables. As
      before the obvious reason is that there are exactly two entries
      where relevant perturbations can occur and that at these entries
      only two resp. one variables really matter for the
      perturbations.\\
\item The most surprising fact is that in 6 variables there are simple
      singularities which are not quasihomogeneous in the strict
      sense, since one of the variables has to be of non-positive
      weight in order to satisfy relative row and column weights.
\end{itemize}

\section{Singularities in $({\mathbb C}^2,0)$ - Fat Points}
\label{FatPoints}
As mentioned at the beginning of the previous section, the usual counting
argument does not provide a bound for the number of rows of a candidate
in dimension 2. Therefore the process of finding candidates has to consider
all possible matrix sizes. Thus we start by considering the $3 \times 2$
matrices, try to find weighted jets which do not allow simple singularities
of this size and mark the remaining cases as candidates. Iterating this, we
increase the matrix size step by step until we reach a matrix size for
which we can prove that there cannot be any simple singularities
(which luckily happens for the $4 \times 3$ matrices).

The following lemma will be the most important tool for proving that
certain candidates are not simple in the case of fat points:

\begin{lem}
\label{dim2l1}
Quasihomogeneous $3 \times 2$ matrices w.r.t. the weights
$$\left( \begin{pmatrix}1 & 1 & 1 \cr 2 & 2 & 2\end{pmatrix}, \begin{pmatrix}1 & 1\end{pmatrix} \right)$$
respectively the weights
$$\left(\begin{pmatrix}1 & 1 & 2 \cr 1 & 1 & 2\end{pmatrix}, \begin{pmatrix}1 & 1\end{pmatrix}\right)$$
cannot define a simple fat point singularity in $({\mathbb C}^2,0)$.
The same statement holds for $4 \times 3$ matrices w.r.t. the weights
$$\left( \begin{pmatrix}1 & 1 & 1 & 1 \cr
                 1 & 1 & 1 & 1 \cr
                 2 & 2 & 2 & 2\end{pmatrix}, \begin{pmatrix}1 & 1\end{pmatrix} \right)$$
respectively
$$\left(\begin{pmatrix}1 & 1 & 1 & 2 \cr
                 1 & 1 & 1 & 2 \cr
                 1 & 1 & 1 & 2\end{pmatrix}, \begin{pmatrix}1 & 1\end{pmatrix}\right).$$
\end{lem}

\begin{proof}
In the first case, the total number of monomials which can appear in the
matrix is $\sum_{ij} r(a,D_{ij})=15$ and the cardinalities of the sets $S_i$
are $S_1=2$, $S_2=9$ and $S_3=4$ which shows that
$$\#S_1 +\#S_2+\#S_3 -2 = 13 < 15 = \sum_{ij} r(a,D_{ij}).$$
Therefore this singularity cannot be simple. \\
In the other cases the inequalities are
$$\#S_1 +\#S_2+\#S_3 -2 = 13 < 14 = \sum_{ij} r(a,D_{ij}),$$
$$\#S_1 +\#S_2+\#S_3 -2 = 25 < 28 = \sum_{ij} r(a,D_{ij}),$$
$$\#S_1 +\#S_2+\#S_3 -2 = 25 < 27 = \sum_{ij} r(a,D_{ij}).$$
\end{proof}

As in the previous section, we now consider the 1-jets:

\begin{lem}
\label{dim2l2}Let $M$ be a $3 \times 2$ matrix with entries in the maximal
ideal of ${\mathbb C}\{x,y\}$. Then $j_1M$ is contact-equivalent to one of
the two jets
$$\begin{pmatrix}x & y & 0 \cr 0 & 0 & y\end{pmatrix}
\text{ and }
\begin{pmatrix}x & y & 0 \cr 0 & x & y\end{pmatrix}$$
or $M$ cannot be simple.
\end{lem}

\begin{proof}
Without loss of generality, we may assume that the matrix of $j_1M$ is
either of the form
$$\begin{pmatrix}x & y & 0 \cr \ast & \ast & \ast\end{pmatrix}
\text{ or }
\begin{pmatrix}x & 0 & 0 \cr \ast & \ast & \ast\end{pmatrix},$$
where the $\ast$ denotes an entry which may be zero or of degree 1.

In the {\bf first} case, the last entry in the bottom row has to be non-zero,
because otherwise $M$ would not be simple by adjacency to a matrix of the
second case of the previous lemma. By exchanging the roles of $x$ and $y$
and permuting the corresponding columns if necessary, we may assume that
the 1-jet is
\begin{eqnarray*}
\begin{pmatrix} x & y & 0 \cr \ast & \ast & y + \alpha x\end{pmatrix} & \sim_C &
  \begin{pmatrix} x & y-\alpha x & 0 \cr \ast & \ast & y\end{pmatrix}\\
 & \sim_C &
  \begin{pmatrix} x & y & 0 \cr \gamma_1 x & \gamma_2 x & y\end{pmatrix} \\
 & \sim_C &
  \begin{pmatrix} x & y & 0 \cr 0 & \gamma_2 x & y\end{pmatrix}\\
\end{eqnarray*}
where $\alpha$,$\gamma_1$ and $\gamma_2 \in {\mathbb C}$. Moreover, we can
assume at this point that $\gamma_2 \in \{0,1\}$ by multiplying the last
row and last colum n by suitable constants if necessary.

If $\gamma_2 = 0$, $M$ has the form of the first matrix in the statement of the
lemma, of the second matrix otherwise.

In the {\bf second} case from the beginning of the proof, we may again
assume by the considerations from above that the 1-jet is
$$\begin{pmatrix} x & 0 & 0 \cr \ast & \ast & \alpha_1 x +\alpha_2 y\end{pmatrix}$$
\begin{entry}
\item[$\alpha_2 \neq 0$:]
By direct computation, we obtain that the 1-jet is contact-equivalent to
$$\begin{pmatrix} x & 0 & 0 \cr 0 & \beta x & y\end{pmatrix}$$
where $\beta \in \{0,1\}$. If $\beta=0$, $M$ cannot be simple by the previous
lemma, otherwise we are in case 1 of the statement of this lemma.

\item[$\alpha_2 = 0$:]
By direct computation, it turns out that this case is up to permutation of
columns identical to case 1.
\end{entry}
\end{proof}

\begin{lem}
Let $M$ be a $4 \times 3$ matrix with entries in the maximal ideal of
${\mathbb C}\{x,y\}$. Then $M$ cannot be simple.
\end{lem}

\begin{proof}
After suitable row and column operations we may assume that $j_1M$ is either
of the form
$$\begin{pmatrix}x & \alpha y & 0 & 0 \cr
           0 & \ast & \ast & \ast \cr
           0 & \ast & \ast & \ast\end{pmatrix}
\text{ or }
\begin{pmatrix}x & \alpha y & 0 & 0 \cr
           y & \ast & \ast & \ast \cr
           0 & \ast & \ast & \ast\end{pmatrix},$$
where $\alpha \in \{0,1\}$.

In the first case, let us consider the following perturbations of $M$ resp.
$j_1M$:
\begin{eqnarray*}
\begin{pmatrix}x & \alpha y & 0 & t \cr
         0 & \ast & \ast & \ast \cr
         0 & \ast & \ast & \ast\end{pmatrix} & \sim_{C \text{ of\;\; 1-jets}} &
\begin{pmatrix}0 & 0 & 0 & t \cr
         0 & \ast & \ast & 0 \cr
         0 & \ast & \ast & 0\end{pmatrix} \\
 & \sim_{C \text{ of\;\; 1-jets}} &
\begin{pmatrix}0 & \ast & \ast \cr
         0 & \ast & \ast\end{pmatrix}
\end{eqnarray*}
which is by lemma \ref{dim2l1} not a 1-jet of the matrix of a simple fat point.

In the other case, we have to consider 2 subcases, namely
\begin{eqnarray*}
j_1M & \sim_C & \begin{pmatrix}x & \alpha y & 0 & 0 \cr
                          y & \ast & \ast & y \cr
                          0 & \ast & \ast & \ast\end{pmatrix}\\
 & \sim_C & \begin{pmatrix}x & \alpha y & 0 & 0 \cr
                      0 & \ast & \ast & y \cr
                      0 & \ast & \ast & \ast\end{pmatrix}\\
\end{eqnarray*}
(which is again in the first case from the beginning of this proof) or
$$j_1M \sim_C \begin{pmatrix}x & \alpha y & 0 & 0 \cr
                        y & \ast     & 0 & x \cr
                        0 & \ast     & \ast & \ast\end{pmatrix}$$
In this last subcase, we see from the third column that either the singularity defined by $M$ cannot be simple by lemma \ref{dim2l1} (if the third column is zero),
or $j_1M$ is again in the first case of the proof by exchanging appropriate columns and rows.
\end{proof}

Therefore we see that a simple fat point singularity in $({\mathbb C}^2,0)$
which is not a complete intersection has to be described by a matrix
with 1-jet
$$\begin{pmatrix}x & y & 0 \cr 0 & \beta x & y\end{pmatrix}$$
where $\beta \in \{0,1\}$. By the same arguments as in the proof of lemma
\ref{dim2l2}, we see that the matrix will then be contact-equivalent to
$$\begin{pmatrix}x & y & 0 \cr 0 & x^k & y\end{pmatrix}.$$

\begin{lem}
\label{lem4l4}
The only fat point singularites in $(\C^2,0)$ which are simple but not complete intersections are listed in the following table:
\begin{center}
\begin{tabular}{|c|c|c|c|}
\hline Type & Equations & $\tau$ &  \struta \\
\hline \strutb $\Xi_k$ & $\begin{pmatrix} x & y & 0 \\ 0 & x^k & y \end{pmatrix}$ &
 $k+3$ & $k\geq 1$  \\
\hline
\end{tabular}
\end{center}
\end{lem}

\begin{proof}
Using $T^1_X$ , as in the last chapter, we can determine the versal family
and calculate the possible deformations of $M$:
$$M'\sim_C\begin{pmatrix}x & y+\beta & \gamma \cr \alpha & x^k+\delta_{k-1}x^{
k-1}+\ldots+\delta_0 & y\end{pmatrix}$$ with
$\alpha,\beta,\gamma,\delta_{k-1},\ldots,\delta_{0}\in\C$.
\\
\begin{multline*}
\begin{pmatrix}x & y+\beta & \gamma \cr \alpha & x^k+\delta_{k-1}x^{k-1}+\ldots+\delta_0 & y\end{pmatrix} \\
 \sim_C \left\{
\begin{array}{ccll}
\left<x^{k+2},y\right> & \text{ for } & \alpha\neq0 & A_{k+1} \\
\left<xy,x^k+y2\right> & \text{ for } & \alpha=0, \gamma\neq0 & F^{k,2
}_{k+1} \\
\left<x,y^2\right> & \text{ for } & \alpha=\gamma=0, \delta_0\neq0 & A_1
\\
\left<x^{k+1},y\right> & \text{ for } & \alpha=\gamma=\delta_0=0, \beta\neq 0 & A_k \\
\begin{pmatrix}x & y & 0 \cr 0 & x^i & y\end{pmatrix} & \text{ for } & \alpha=\gamma=\delta_0=\gamma=0 & \Xi_i, i\leq k \\
\end{array}
\right.
\end{multline*}
Since all (finitely many types of) singularities appearing here are simple,
we have proved that $\Xi_k$ is simple.
\end{proof}

\clearpage
\section{Adjacencies}

When proving a classification of simple isolated singularities, a large number
of adjacencies is usually determined explicitly as a byproduct or follows
from these by transitivity (see \cite{Arn}, \cite{Giu}, \cite{FK1}).
In particular, knowing the complete list of adjacencies for a given series
(resp. at least knowing into which series it can deform into) is vital
for proving that all singularities in this series are simple.
  
On the other hand, when proving that a singularity, which is not part of
a series, is indeed simple, it is often sufficient to know by semicontinuity
of certain invariants (e.g. $\delta$, $\mu$) that adjacencies to non-simple
singularities cannot occur.\\

Due to the large number of simple singularities in 4 and more
variables, we only concentrate on 2 and 3 variables for computing
adjacencies where the case of 3 variables may be reagarded as a model
for the cases with more than 3 variables. In order to determine the
complete list of adjacencies for isolated simple fat point and space
curve singularities singularities, we will now start by summarizing
what is known:

\noindent
While the list of adjacencies from Arnold's classification \cite{Arn}
can easily be checked to be complete, there are known gaps in Giusti's
list, as was first observed when additional adjacencies were found for
space curves singularities by Goryunov \cite{Gor}. These gaps were closed
in \cite{S-V} for the adjacencies among the isolated simple complete
intersection space curve singularities. The case of adjacencies to plane
curve singularities has later been treated by computer algebra methods
in \cite{FK2}.
\begin{center}
\includegraphics[width=\textwidth]{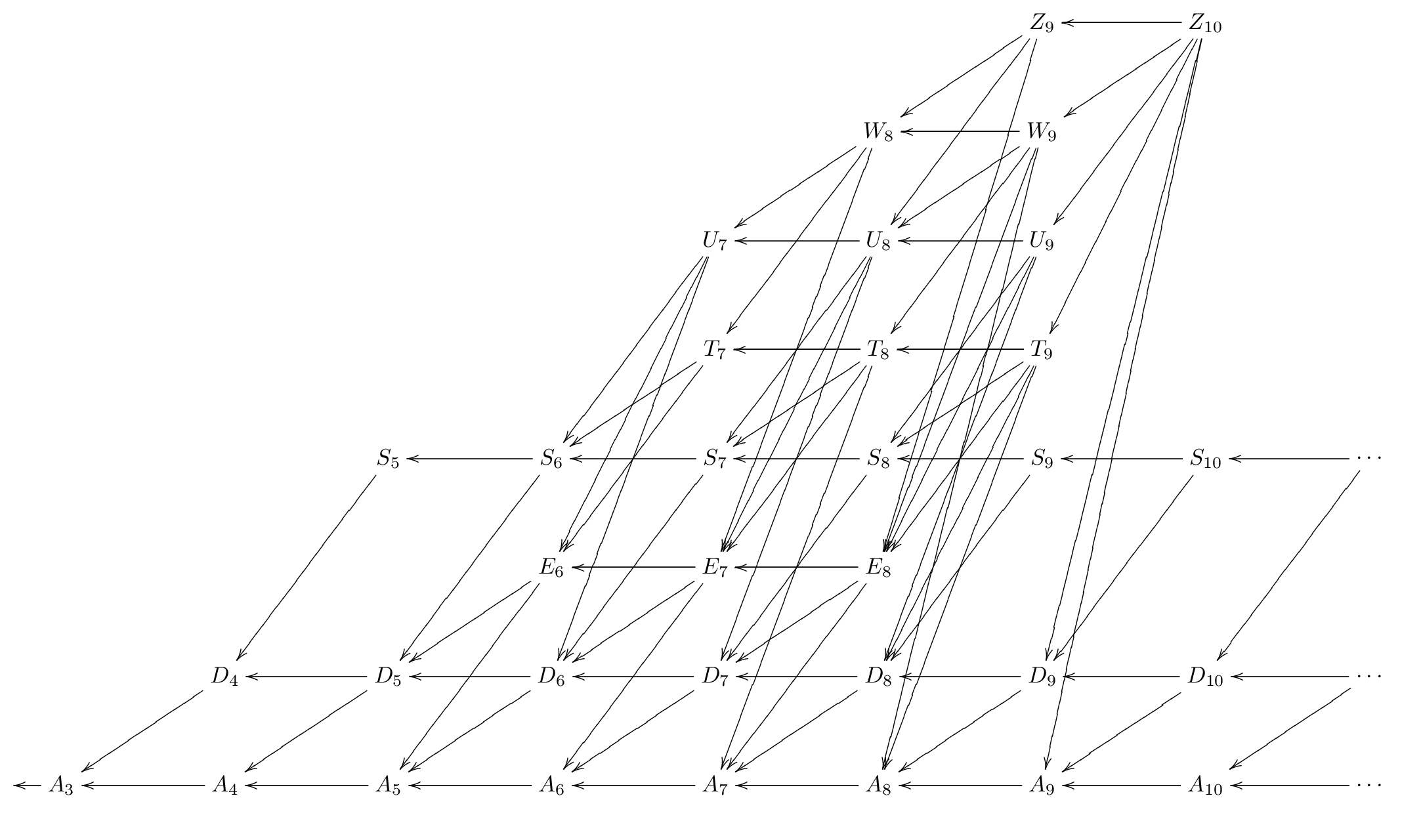}\\
Diagram 1: Adjacencies of simple ICIS space curve singularities
due to
\cite{Arn},\cite{Giu},\cite{Gor},\cite{S-V},\cite{FK2}
\end{center}

In the case of simple fat point singularities, the adjacencies of the series $\Xi_k$ already follow immediately from the proof of lemma \ref{lem4l4}. So the
remaining considerations in this case concern Giusti's list. Here it turns out that the adjacency table given in \cite{Giu} is not just incomplete, but even
contains non-existent adjacencies. Therefore we state corrected versions of the adjacency tables for the series $F$, $H$ and $I$ here\footnote{The adjacency lists
for the singularities of types $G_5$ and $G_7$ did not contain any mistakes.}:
\\
\begin{figure}
\begin{flushright}
\includegraphics[angle=90,width=0.75\textwidth]{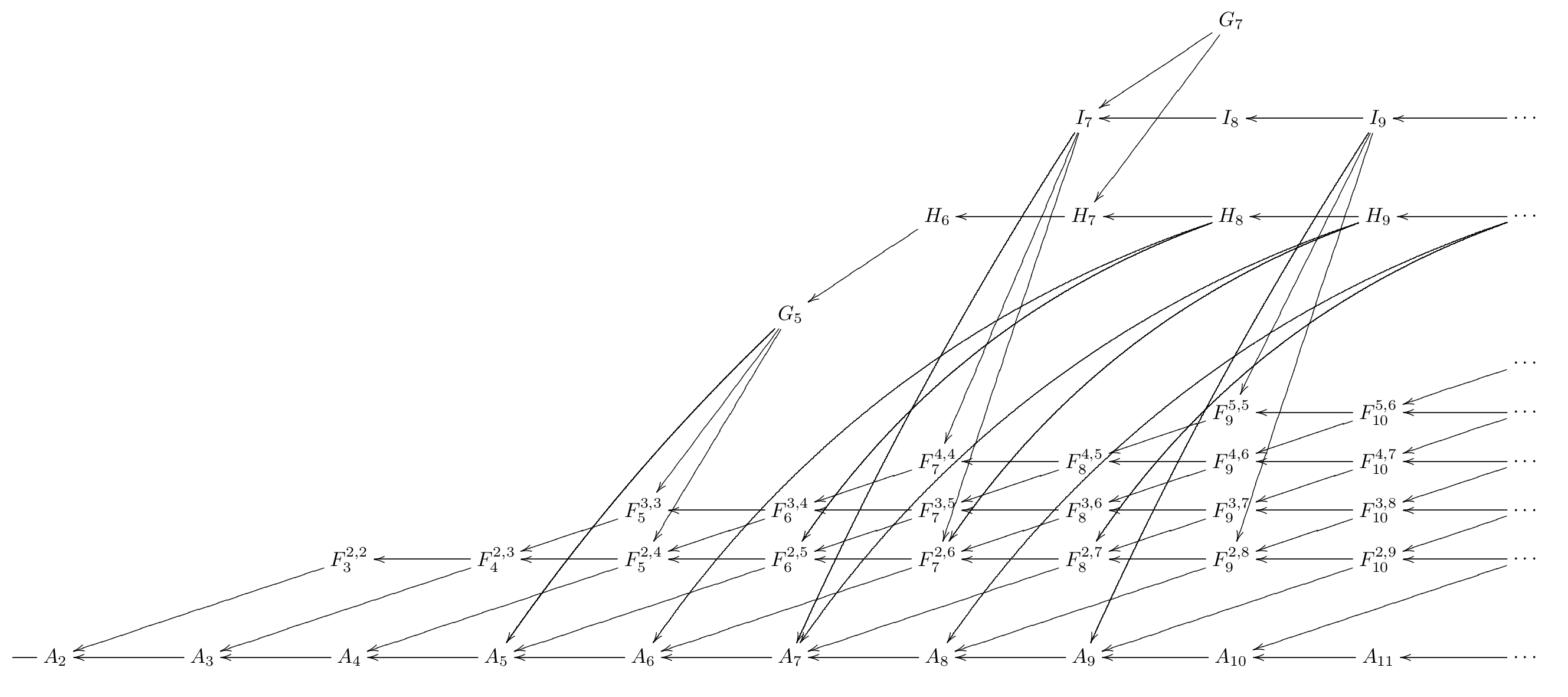}
\end{flushright}
\begin{center}
Diagram 4: Adjacencies of simple ICIS fat point singularities
from Giusti's classification \cite{Giu} corrected according to the adjacency lists
\end{center}
\end{figure}

\vspace{6pt}\noindent Adjacencies for singularities of type $F^{n,p}_{n+p-1}$:
{\small 
\begin{multline*}
\left<\alpha x+\beta y+xy,a_1x+\ldots+x^n+b_1y+\ldots+y^p\right>, \\
 \sim_C \left\{
\begin{array}{ccll}
\left<x,y^p\right> & \text{ for } & \alpha\neq 0 & A_{p-1} \\
\left<y,x^n\right> & \text{ for } & \alpha=0, \beta\neq 0 & A_{n-1} \\
\left<y^{p+1},x\right> & \text{ for } & \alpha=\beta=0, a_1\neq 0 & A_p \\
\left<x^{n+1},y\right> & \text{ for } & \alpha=\beta=a_1=0, b_1\neq 0 & A_n \\
\left<xy,x^i+y^j\right> & \text{ for } & \alpha=\beta=a_1=b_1=0 & F^{i,j}_{i+j-1}, i\leq n, j\leq p
\end{array}
\right.
\end{multline*}}

\vspace{20pt}\noindent
Adjacencies for singularities of type $H_{n+3}$:
{\small
\begin{multline*}
\left<x^2+\alpha x+a_1y+\ldots+a_{n-1}y^{n-1}+y^n, xy^2+\beta x+\gamma y+bxy+cy^2+dy^3\right>\\
 \sim_C \left\{
\begin{array}{ccll}
\left<x,y^{n+2}\right> & \text{ for } & \alpha\neq 0 & A_{n+1} \\
\left<y,x^5\right> & \text{ for } & \alpha=0, a_1\neq 0 & A_4 \\
\left<y^n,x\right> & \text{ for } & \alpha=a_1=0, \beta\neq 0 & A_{n-1} \\
\left<x^2,y\right> & \text{ for } & \alpha=a_1=\beta=0, \gamma\neq 0 & A_1 \\
\left<x^2+y^2,xy\right> & \text{ for } & b\neq 0, a_2\neq c^2 & F^{2,2}_3 \\
\left<x^2+y^3,xy\right> & \text{ for } & b\neq 0, a_2=c^2, a_3\neq cd & F^{2,3}_4 \\
\left<x^2+y^4,xy\right> & \text{ for } & b\neq 0, a_2=c^2, a_3=cd, a_4\neq d^2 & F^{2,4}_5 \\
\left<x^2+y^i,xy\right> & \text{ for } & b\neq 0, a_2=c^2, a_3=cd, a_4=d^2,\\
 & &  a_5=\ldots=a_{i-1}=0,a_i\neq 0, i\leq n & F^{2,i}_{i+1} \\
\left<x^2+y^2,xy\right> & \text{ for } & b=0, a_2\neq 0, c\neq 0 & F^{2,2}_3 \\
\left<x^3+y^3,xy\right> & \text{ for } & b=0, a_2\neq 0, c=0 & F^{3,3}_5 \\
\left<x^2+y^2,xy\right> & \text{ for } & b=a_2=0, c\neq 0 & F^{2,2}_3\\
\left<x^2,y^3\right> & \text{ for } & b=a_2=c=0, d\neq 0 & G_5\\
\left<x^2+y^{n-1},xy^2\right> & \text{ for } & b=a_2=\ldots=a_{n-2}=c=d=0, & H_{n+2}\\
 & & a_{n-1}\neq 0 & \\
 \end{array}
\right.
\end{multline*}}

\vspace{20pt}\noindent
Adjacencies for singularities of type $I_{2p-1}$:
{\small
\begin{multline*}
\left<x^2+y^3+ay^2+\alpha x+\beta y,\right.\\
\left.\gamma x+\delta y+b_1xy+\ldots+b_{p-2}xy^{p-2}+c_2y^2+\ldots+c_{p-1}y^{p-1}+y^p\right> \\
 \sim_C \left\{
\begin{array}{ccll}
\left<x,y^{p+1}\right> & \text{ for } & \alpha\neq 0 & A_p \\
\left<y,x^{2p}\right> & \text{ for } & \alpha=0, \beta\neq 0 & A_{2p-1} \\
\left<y^{2p},x\right> & \text{ for } & \alpha=\beta=0, \gamma\neq 0 & A_{2p-1} \\
\left<x^2,y\right> & \text{ for } & \alpha=\beta=\gamma=0, \delta\neq 0 & A_1 \\
\left<x^2+y^i,xy\right> & \text{ for } & b_1\neq 0, i\leq 2p-2 & F^{2,i}_{i+1} \\
\left<x^i+y^j,xy\right> & \text{ for } & b_1=0, a\neq 0, i,j\leq p & F^{i,j}_{i+j-1} \\
\left<x^2+y^2,xy\right> & \text{ for } & b_1=a=0, c_2\neq 0 & F^{2,2}_3 \\
\left<x^2,y^3\right> & \text{ for } & b_1=a=c_2=0, c_3\neq 0 & G_5 \\
\left<x^2+y^3,xy^2\right> & \text{ for } & b_1=a=c_2=c_3=0, b_2\neq 0 & H_6 \\
\left<x^2+y^3,xy^{p-2}\right> &  & \text{otherwise }  & I_{2p-2}
\end{array}
\right.
\end{multline*}}

\vspace{30pt}\noindent
Adjacencies for singularities of type $I_{2q+2}$:
\begin{multline*}
\left<x^2+ay^2+y^3+\alpha x+\beta y,\right.\\
\left.\gamma x+\delta y+b_1xy+\ldots+b_{q-1}xy^{q-1}+xy^q+c_2y^2+\ldots+c_{q+1}y^{q+1}\right>  \\
 \sim_C \left\{
\begin{array}{ccll}
\left<x,y^{q+3}\right> & \text{ for } & \alpha\neq 0 & A_{q+2} \\
\left<y,x^{2q+2}\right> & \text{ for } & \alpha=0, \beta\neq 0 & A_{2q+1} \\
\left<y^{2q+2},x\right> & \text{ for } & \alpha=\beta=0, \gamma\neq 0 & A_{2q+1} \\
\left<x^2,y\right> & \text{ for } & \alpha=\beta=\gamma=0, \delta\neq 0 & A_1 \\
\left<x^2+y^i,xy\right> & \text{ for } & b_1\neq 0, i\leq 2q & F^{2,i}_{i+1} \\
\left<x^i+y^j,xy\right> & \text{ for } & b_1=0, a\neq 0, i,j\leq q+1 & F^{i,j}_{i+j-1} \\
\left<x^2+y^2,xy\right> & \text{ for } & b_1=a=0, c_2\neq 0 & F^{2,2}_3 \\
\left<x^2,y^3\right> & \text{ for } & b_1=a=c_2=0, c_3\neq 0 & G_5 \\
\left<x^2+y^3,xy^2\right> & \text{ for } & b_1=a=c_2=c_3=0, b_2\neq 0 & H_6 \\
\left<x^2+y^3,y^{q+1}\right> &  & \text{otherwise }  & I_{2q+1}
\end{array}
\right.
\end{multline*}

\clearpage
\vspace{20pt}
\begin{center}
\includegraphics[width=\textwidth]{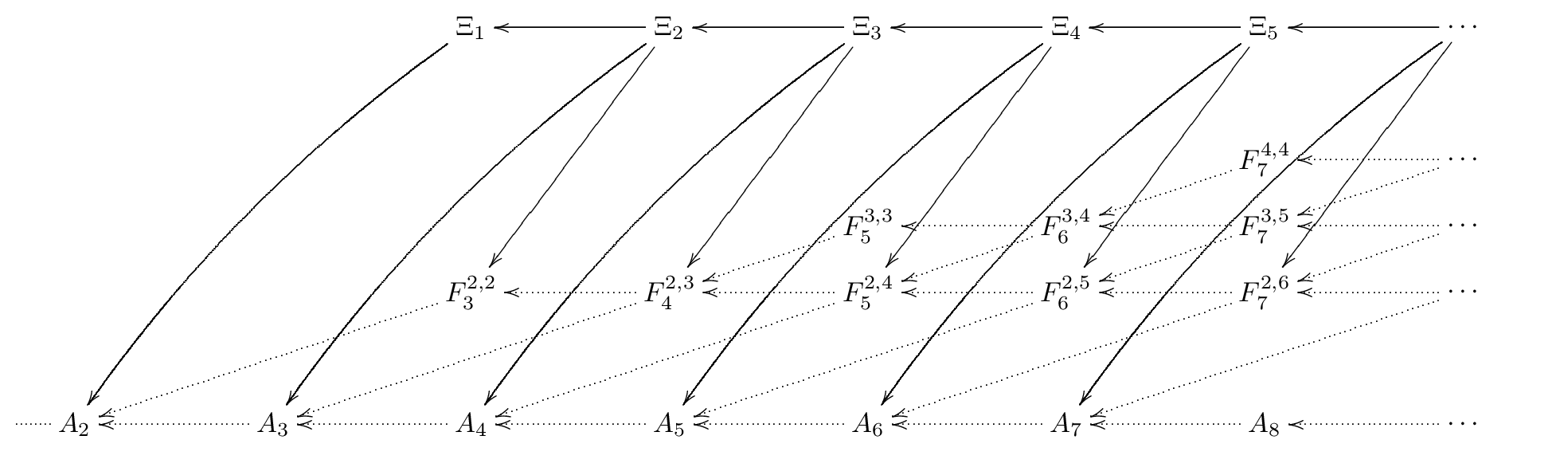}
Diagram 5: Adjacencies of simple non-complete-intersection fat point singularities. As singularities of type $\Xi_k$ cannot deform into singularities of types $G$,
$H$ and $I$, only singularity types $F$ and $A$ appear in this diagram, adjacencies among those are drawn as dotted lines. \vspace{10pt}
\end{center}

The last item that remains to be determined is the complete list of
adjacencies for simple isolated space curve singularities which are not
complete intersections\footnote{The classification of these singularities can
be found in \cite{FK1}}. Due to the large number of simple space curve
singularities, this turns out to be a rather lengthy calculation which
should begin at the singularities with the lowest values for the invariants
$\tau$ and $\delta$ and always involves the same steps for each given simple
singularity. Therefore we only sketch the basic concept of this calculation,
the detailed results for each singularity and each step can be found at
\begin{center}
http://www.mathematik.uni-kl.de/$\sim$anne/adjCMcod2.html
\end{center}

\noindent
\begin{entry}
\item[Step 1:] Determining candidates for new adjacencies\\
For a given singularity, we first consider all singularities, whose
Tjurina number is smaller than the one of the given singularity, as possible
target of an adjacency. From this list, we exclude, of course, those
singularities for which it has been shown in the proof of simplicity that
the adjacency does not exist and mark the ones which have been shown to exist.

Using the semicontinuity of the invariants $\delta$ and $\mu$ and the fact
that $\delta$ has to be constant in a $\mu$-constant family (cf. \cite{B-G}),
we can then exclude some more singularities from the list.
In the case of the singularity $S_6^*$, for example, the invariants are
$\tau = 7$, $\mu = 6$ and $\delta =4$ and thus adjacencies to the singularities
$A_6$ ($\mu = 6$, $\delta = 3$) and $E_6$ ($\mu = 6$, $\delta$ = 3) are
excluded by the latter condition.\\
\\
\item[Step 2:] Finding adjacencies for which $\tau$ drops exactly by one\\
For finding these adjacencies, we have to study the structure of the base
space of the versal family of the given singularity more closely. To this
end, let $t_1,\dots,t_{\tau}$ denote the parameters of the versal family
and consider the relative $T^1$ of the family as a
$\C[t_1,\dots,t_{\tau}]$-module. The flattening stratification\footnote{From
the computational point of view, determining the flattening stratification
means computing a presentation matrix of the module, which is described in
\cite{FK2}, and then forming the Fitting ideals. Since computing minors of
matrices with polynomial entries is quite expensive, we first only consider
the stratum on which $\tau$ drops by exactly one which is given by the
vanishing of the 2-minors and non-vanishing of the 1-minors of the presentation
matrix. It turns out during the computation that we can actually avoid
studying the stratum on which $\tau$ drops by 2 in all cases.} of this
module determines the strata in the base space of the versal family on which
the Tjurina number is constant. \\
Since all simple singularities from our list
are quasihomogeneous, we know that we also have an Euler relation for the
given singularity and thus the ($\tau$-1)-th Fitting ideal of
the relative $T^1$ is the maximal ideal at the origin. The stratum on
which $\tau$ drops exactly by one is determined by
${\rm Fitt}_{(\tau -2)}(T^1_{\it rel})$
ideal on the complement of ${\rm Fitt}_{(\tau -1)}(T^1_{\it rel})$.\\
On each primary component of ${\rm Fitt}_{(\tau -2)}(T^1_{\it rel})$, we
have exactly one type of singularity. So, it now suffices to determine a
primary decomposition of this ideal and determine the type of singularity
on each component. The types of singularities appearing there are exactly
those of the correct Tjurina number $\tau$ to which an adjacency exists.
This allows us to exclude resp. to mark further singularities in our list
of candidates (of course also using the transitivity of adjacencies again).\\
\\
During the computations for all simple space curve singularities, it turned
out that after this step no unmarked candidates were left in all cases
except for the singularities $E_k \vee L$, $k \in \{6,7,8\}$.\\
\\
\item[Step 3:] Treating the remaining cases\\
For the singularities $E_k \vee L$, $k \in \{6,7,8\}$, there is only one
remaining candidate for an adjacency, namely $E_k \vee L \longrightarrow A_k$.
To exclude this adjacency, we first observe that a versal family of an
$E_k \vee L$ singularity is of the form
$$\begin{pmatrix}
z & \alpha & \beta \cr
t_1 & x & y
\end{pmatrix},$$
where $\alpha$ and $\beta$ are suitable polynomials in $x$, $y$ and the
parameters
$t_2,\dots,t_{\tau}$. Moreover, it is easy to see that a plane curve
singularity which is not an $A_1$ can only occur if $t_1 \neq 0$. But for
$t_1 \neq 0$ and $t_2=\dots=t_{\tau}=0$, we obtain a $U_{k+1}$ singularity.
Moreover, any singularity with $t_1 \neq 0$ in the versal family of
$E_k \vee L$ also appears in the versal family of $U_{k+1}$ as one can check
by direct computation. Now we already know from \cite{FK2} that $U_{k+1}$
cannot deform into $A_k$. This excludes the last remaining candidates.
\end{entry}
\begin{figure}[p]
\begin{flushright}
\includegraphics[angle=90,height=0.9\textheight]{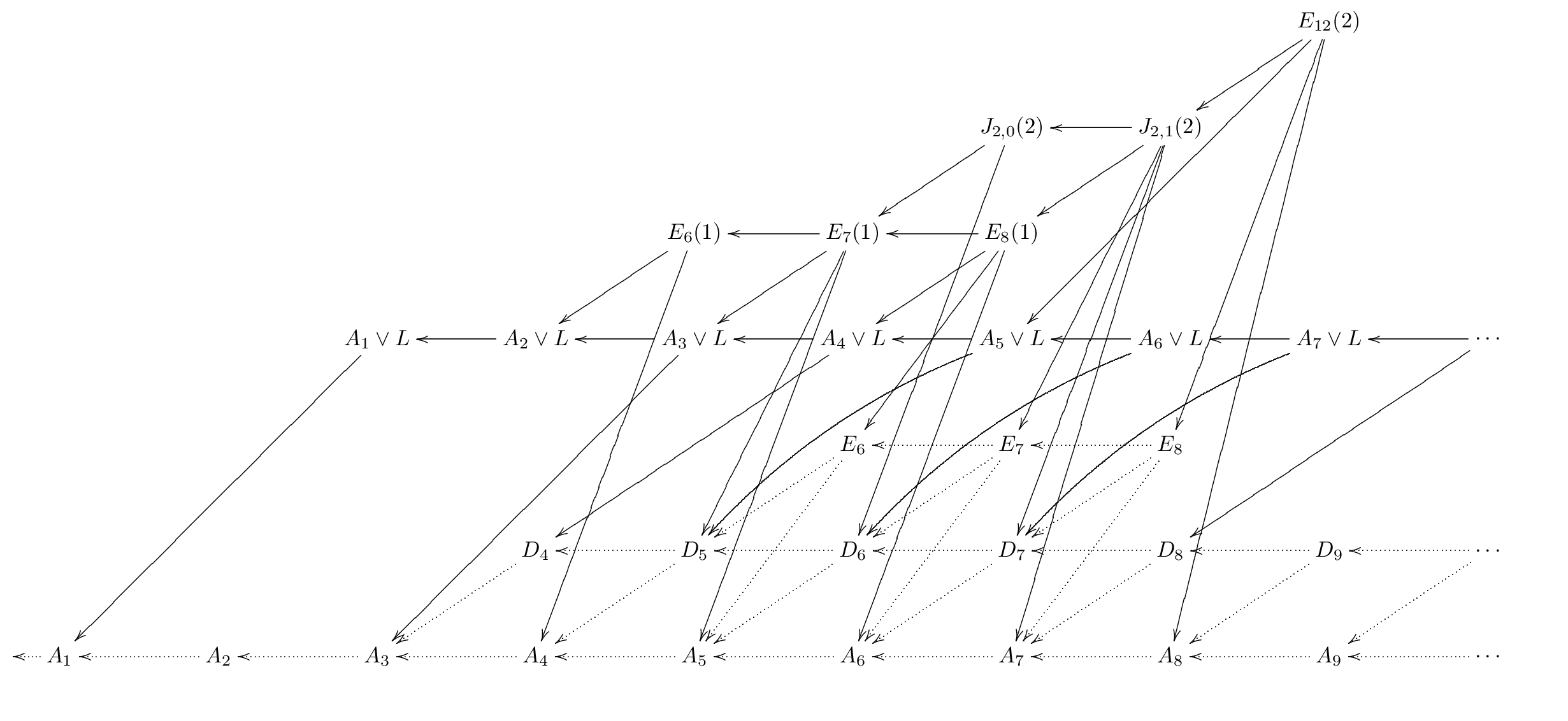}
\end{flushright}\begin{center}Diagram 6: Adjacencies of simple space curve
singularities of multiplicity 3
(Adjacencies of plane curves from Arnold's list are shown as
dotted lines.)
\end{center}
\end{figure}
\begin{figure}[p]
\begin{flushright}
\includegraphics[angle=90,width=\textwidth]{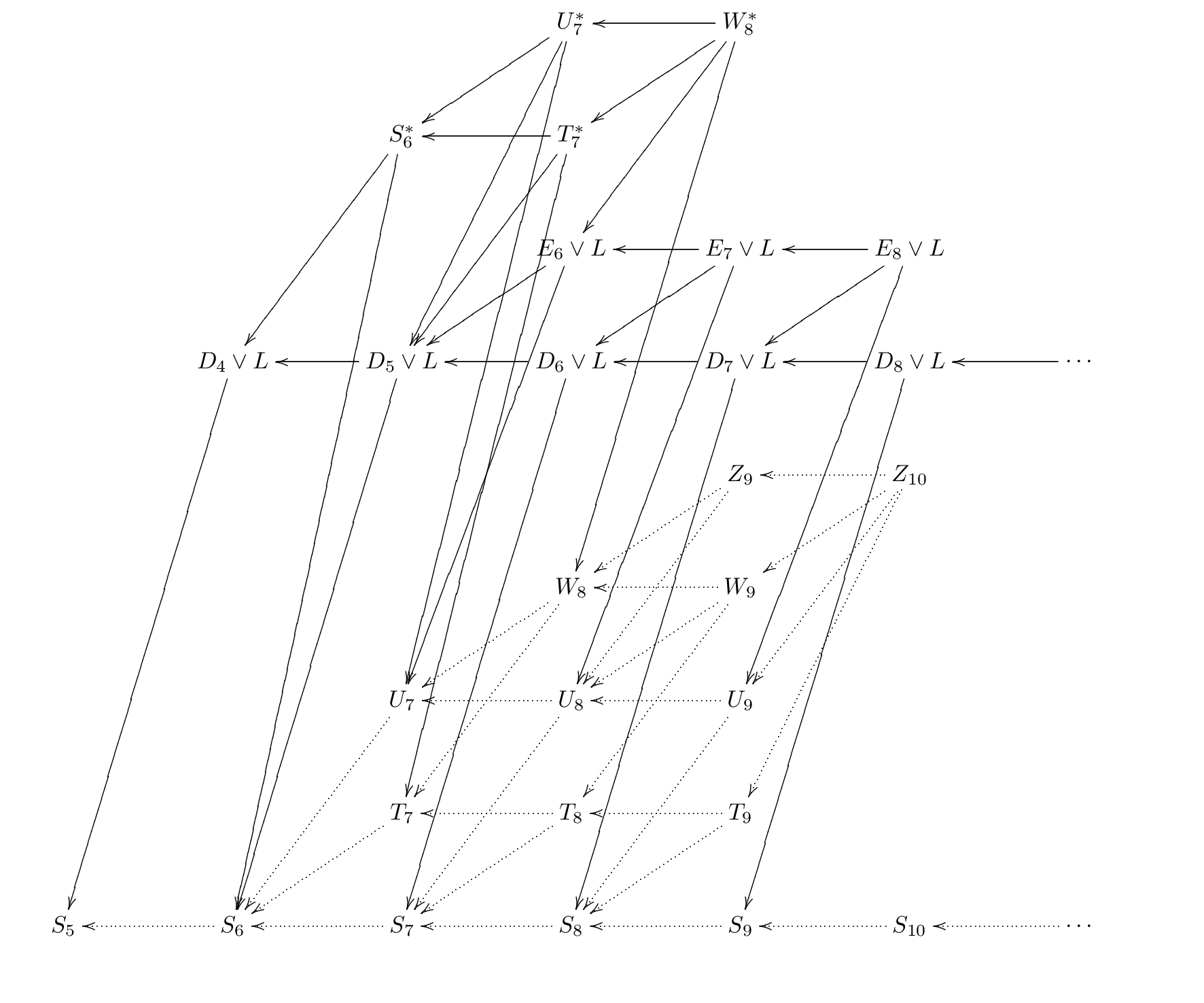}
\end{flushright} \begin{center}Diagram 7: Adjacencies among simple space curve
singularities of multiplicity
4 to space curve singularities of multiplicity 4 (Adjacencies from
\cite{Giu}, \cite{S-V}, \cite{FK2} are shown as dotted lines.)
\end{center}
\end{figure}
\begin{figure}[p]
\begin{flushright}
\includegraphics[angle=90,height=0.85\textheight]{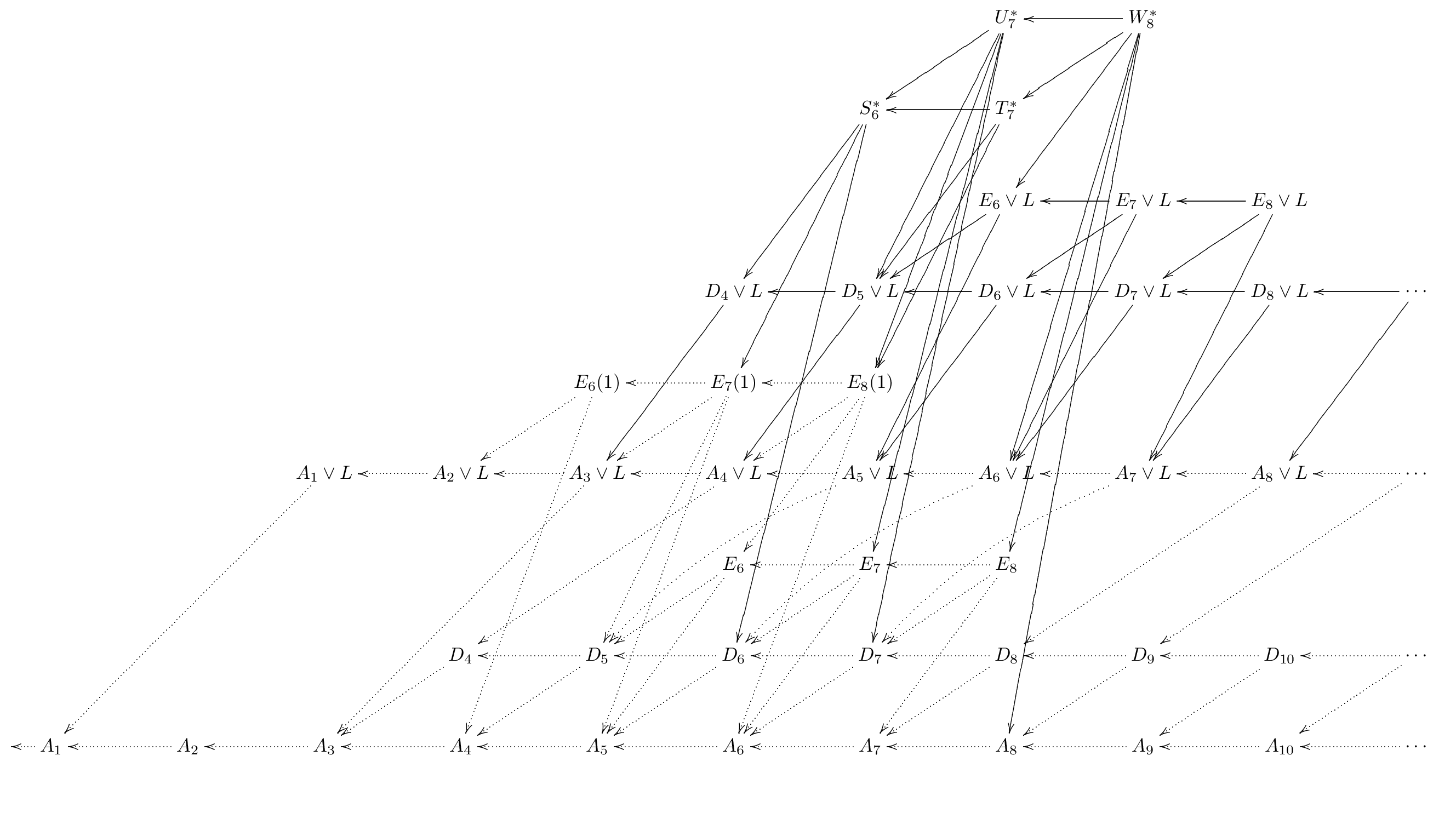}
\end{flushright} \begin{center}Diagram 8: Adjacencies of simple space curve
singularities of multiplicity
4 into singularities of multiplicity at most 3 (Adjacencies from
\cite{Arn}, \cite{FK2} are shown as dotted lines. Note that
further adjacencies to singularities of multiplicity 3 do exist as
combination of an adjacency of the previous diagram and diagram 1
due to transitivity )
\end{center}
\end{figure}

\clearpage
\section{Appendix} For readers' convenience, Arnold's list of simple hypersurface singularities from \cite{Arn} and
the complete lists of previously known simple Cohen-Macaulay
codimension 2 singularities (cf. \cite{Giu} and \cite{FK1}) are listed
in this Appendix.

\vspace{0,3cm}

\begin{center}
\begin{tabular}{|c|c|c|c|c|c|}
\hline Type & Presentation Matrix  & $\mu$ & $\tau$ & $\delta$ & \struta
\\ \hline
$A_k$ & $x^{k+1}+y^2$  & $k$ & $k$ & $\left\lfloor\frac{k}{2}\right\rfloor$   &  $k\geq 1$\struta \\
\hline
$D_k$ & $x^2y+y^{k-1}$ & $k$ & $k$ & $\left\lfloor\frac{k+2}{2}\right\rfloor$ &  $k\geq 4$ \struta \\
\hline
$E_6$ & $x^3+y^4$ & 6 & 6 & 3  & \struta \\
\hline
$E_7$ & $x^3+xy^2$ & 7 & 7 & 4  & \struta \\
\hline
$E_8$ & $x^3+y^5$ & 8 & 8 & 4  & \struta \\
\hline
\end{tabular}
\\
\vspace{10pt}
Table 1: The simple hypersurface singularities
\vspace{30pt}
\\
\begin{tabular}{|c|c|c|c|c|}
\hline Type & Presentation Matrix  & $\mu$ & $\tau$ & \struta
\\ \hline
$A_k$ & $\left<y,x^{k+1}\right>$  & $k$ & $k$ & $k\geq 1$ \struta\\
\hline
$F_{q+r-1}^{q,r}$ & $\left<xy,x^q+y^r\right>$ & $q+r-1$ & $q+r$   &  $q,r\geq 2$ \struta\\
\hline
$G_5$ & $\left<x^2,y^3\right>$ & 5 & 7 & \struta \\
$G_7$ & $\left<x^2,y^4\right>$ & 7 & 10 & \struta \\
\hline $H_{q+3}$ & $\left<x^2+y^{q},xy^2\right>$ & $q+3$ & $q+5$ &
$q\geq 3$\struta
\\
\hline $I_{2q-1}$ & $\left<x^2+y^3,y^q\right>$ & $2q-1$ & $2q+1$ &
$q\geq 4$\struta
\\
$I_{2r+2}$ & $\left<x^2+y^3,xy^r\right>$ & $2r+2$ & $2r+4$ & $r\geq 3$\struta \\
\hline
\strutb $\Xi_k$ & $\begin{pmatrix} x & y & 0 \\ 0 & x^k & y \end{pmatrix}$ &
& $k+3$ & $k\geq 1$  \\
\hline
\end{tabular}
\\
\vspace{10pt}
Table 2: The simple fat point singularities in $(\C^2,0)$
\clearpage

\begin{tabular}{|c|c|c|c|c|c|}
\hline Type & Presentation Matrix  & $\mu$ & $\tau$ & $\delta$ & \struta
\\ \hline
$S_{n+3}$ & $(x^2+y^2+z^{n},yz)$ & $n+3$ & $n+3$ & $\left\lfloor\frac{n+6}{2}\right\rfloor$   &  $n\geq 2$ \struta\\
\hline
$T_7$ & $(x^2+y^3+z^3,yz)$ & 7 & 7 & 4   & \struta \\
\hline
$T_8$ & $(x^2+y^3+z^4,yz)$ & 8 & 8 &  5  & \struta \\
\hline
$T_9$ & $(x^2+y^3+z^5,yz)$ & 9 & 9 & 5   & \struta \\
\hline
$U_7$ & $(x^2+yz,xy+z^3)$ & 7 & 7 & 4   & \struta \\
\hline
$U_8$ & $(x^2+yz+z^3,xy)$ & 8 & 8 & 5   & \struta \\
\hline
$U_9$ & $(x^2+yz,xy+z^4)$ & 9 & 9 & 5   & \struta \\
\hline
$W_8$ & $(x^2+z^3,y^2+xz)$ & 8 & 8 & 4  & \struta \\
\hline
$W_9$ & $(x^2+yz^2,y^2+xz)$ & 9 & 9 & 5   & \struta \\
\hline
$Z_9$ & $(x^2+z^3,y^2+z^3)$ & 9 & 9 & 5   & \struta \\
\hline
$Z_{10}$ & $(x^2+yz^2,y^2+z^3)$ & 10 & 10 & 5   & \struta \\
\hline
\end{tabular}

\vspace{10pt}
Table 3: Simple space curve singularities in  $(\C^3,0)$, Part 1: ICIS
\vspace{20pt}

\begin{tabular}{|c|c|c|c|c|c|}
\hline
\struta Type & Presentation Matrix  & $\mu$ & $\tau$ & $\delta$ &
\\ \hline
\strutb \parbox{1.5cm}{\begin{center}$A_{k-3}\vee L$\\ $k\geq 4$\end{center}} &
$\begin{pmatrix} z & y & x^{k-3} \\ 0 & x & y \end{pmatrix}$ & $k-2$ & $k-1$ & \parbox{1cm}{\begin{center}$\frac{k}{2}$\\ $\frac{k-1}{2}$\end{center}} & \parbox{1cm}{\begin{center}$k$ even\\ $k$ odd\end{center}} \\
\hline
\strutb $E_6(1)$ &  $\begin{pmatrix} z & y & x^2 \\ x & z & y \end{pmatrix}$  & 4 & 5 & 2 &  \\
\hline
\strutb $E_7(1)$ & $\begin{pmatrix} z+x^2 & y & x \\ 0 & z & y \end{pmatrix}$ & 5 & 6 & 3 &  \\
\hline
\strutb $E_8(1)$ & $\begin{pmatrix} z & y & x^3 \\ x & z & y \end{pmatrix}$ & 6 & 7 & 3 &  \\
\hline
\strutb \parbox{1.5cm}{\begin{center}$J_{2,k}(2)$ \\ $k\in\{0,1\}$\end{center}} &
$\begin{pmatrix} z+x^2 & y & x^{k+2} \\ 0 & z & y \end{pmatrix}$ & \parbox{1cm}{\begin{center}6 \\ 7 \end{center} } & \parbox{1cm}{\begin{center}7 \\ 8 \end{center} } & 4  & \parbox{1cm}{\begin{center} $k=0$ \\ $k=1$\end{center} } \\
\hline
\strutb$E_{12}(2)$ & $\begin{pmatrix} z & y & x^3 \\ x^2 & z & y \end{pmatrix}$ & 8 & 9 & 4  & \\
\hline
\strutb \parbox{1.5cm}{\begin{center}$D_{k+4}\vee L$ \\ $k\geq 0$ \end{center} } & $\begin{pmatrix} z & 0
& x^{k+2}-y^2 \\ 0 & x & y
\end{pmatrix}$
 &  $k+5$ & $k+6$ & \parbox{1cm}{\begin{center}$\frac{k+8}{2}$ \\ $\frac{k+7}{2}$ \end{center} } & \parbox{1.5cm}{\begin{center}$k$ even \\ $k$ odd \end{center} }  \\
\hline
\strutb $E_{6}\vee L$ & $\begin{pmatrix} z & -y^2 & -x^3 \\ 0 & x & y \end{pmatrix}$ &  7 & 8 & 4 &  \\
\hline
\strutb $E_{7}\vee L$ & $\begin{pmatrix} z & x^3-y^2 & 0 \\ 0 & x & y \end{pmatrix}$ &  8 & 9 & 5 &  \\
\hline
\strutb $E_{8}\vee L$ & $\begin{pmatrix} z & -y^2 & -x^4 \\ 0 & x & y \end{pmatrix}$ &  9 & 10 & 5 &  \\
\hline
\strutb $S_6^*$ & $\begin{pmatrix} z & x & y \\ 0 & y & x^2-z^2 \end{pmatrix}$ & 6 & 7 & 4 &  \\
\hline
\strutb $T_7^*$ & $\begin{pmatrix} z & x & y \\ 0 & y & x^2-z^3 \end{pmatrix}$ &  7 & 8 & 4 &  \\
\hline
\strutb $U_7^*$ & $\begin{pmatrix} z & xy & x^2 \\ x & z & y \end{pmatrix}$ &  7 & 8 & 4 &  \\
\hline
\strutb $W_8^*$ & $\begin{pmatrix} z & y^2 & x^2 \\ x & z & y \end{pmatrix}$ &  8 & 9 & 4 &  \\
\hline
\end{tabular}

 \vspace{10pt}
Table 4: Simple space curve singularities in  $(\C^3,0)$, Part 2: Non-complete-intersections
\vspace{20pt}
\end{center}


\end{document}